\documentclass[11pt]{article}
\usepackage[utf8]{inputenc}
\usepackage[T1]{fontenc}
\usepackage[frenchb,english]{babel} 
\usepackage{textcomp}
\usepackage{amsmath,amssymb}
\usepackage{amsthm}
                  
\usepackage{lmodern}

\usepackage[a4paper,margin=2cm]{geometry}

\usepackage{graphicx}             
\usepackage{xcolor}               
\usepackage{microtype,stmaryrd}          

\usepackage{hyperref}
\hypersetup{pdfstartview=XYZ}

\usepackage{bbm}
\usepackage{mathrsfs}
\usepackage{subfig}

\def\build#1_#2^#3{\mathrel{
\mathop{\kern 0pt#1}\limits_{#2}^{#3}}}
\def\llbracket{[\hspace{-.10em} [ }
\def\rrbracket{ ] \hspace{-.10em}]}

\newtheorem{theorem}{Theorem}
\newtheorem{proposition}[theorem]{Proposition}
\newtheorem{definition}[theorem]{Definition}
\newtheorem{lemma}[theorem]{Lemma}

\newtheorem{corollary}[theorem]{Corollary}

\def\w{\mathrm{w}}
\def\t{\mathcal{T}}

\def\g{\mathcal{G}}
\def\W{\mathcal{W}}
\def\S{\mathcal{S}}

\def\N{\mathbb{N}}

\def\D{\mathbb{D}}
\def\P{\mathbb{P}}
\def\E{\mathbb{E}}
\def\R{\mathbb{R}}
\def\z{\mathcal{Z}}
\def\y{\mathcal{Y}}

\def\cc{\mathcal{C}}

\def\ve{{\varepsilon}}
\def\la{\longrightarrow}
\def\ov{\overline}
\def\un{\underline}
\def\dd{\mathrm{d}}
\def\wh{\widehat}
\def\wt{\widetilde}

\def\vol{\mathrm{vol}}
\def\ab{{a^{\bullet}}}

\def\bn{\mathbf{n}}

\title{Growth-fragmentation processes in Brownian motion \\
indexed by the Brownian tree\footnote{Supported by the ERC Advanced Grant 740943 {\sc GeoBrown}}}

\author{Jean-Fran\c cois Le Gall, Armand Riera}

\date{\small Universit\'e Paris-Sud}
\begin{document}
\maketitle

\begin{abstract}
We consider the model of Brownian motion indexed by the Brownian tree. For every $r\geq 0$ and every connected component
of the set of points where Brownian motion is greater than $r$, we define the boundary size of this component, and we then show
that the collection of these boundary sizes evolves when $r$ varies like a well-identified growth-fragmentation process.
We then prove that the same growth-fragmentation process appears when slicing a Brownian disk at height $r$ and 
considering the perimeters of the resulting connected components.
\end{abstract}

{\small \tableofcontents}

\section{Introduction}

The main goal of the present work is to prove that the collection of boundary sizes
of excursions of Brownian motion
indexed by the Brownian tree above a fixed level evolves according to
a well-identified growth-fragmentation process when the level increases. Because of the close connections
between Brownian motion
indexed by the Brownian tree and the Brownian map or the Brownian disk, this result also implies that
the collection of boundary sizes of the connected components of the set of points of
a Brownian disk whose distance from the boundary is greater than $r$ evolves 
according to the same growth-fragmentation process. The latter fact may be viewed as
a continuous analog of a recent result of Bertoin, Curien and Kortchemski \cite{BCK} identifying the 
growth-fragmentation process arising as the scaling limit for the collection of 
lengths of cycles obtained by slicing random Boltzmann triangulations
with a boundary at a given height, when the size
of the boundary grows to infinity. In fact, the growth-fragmentation process of \cite{BCK} 
is the same as in our main results, and this strongly suggests that the results of \cite{BCK} could be
extended to more general planar maps with a boundary (see also \cite{BBCK}
for related results).

In order to give a more precise description of our main results, we first need to recall
the notion of Brownian motion indexed by the Brownian tree. The Brownian tree of interest here
is a variant of Aldous' continuum random tree, which is also called the CRT. This tree 
is conveniently defined as the tree $\t_\zeta$ coded by a Brownian excursion $(\zeta_s)_{0\leq s\leq \sigma}$
under the $\sigma$-finite It\^o measure of positive excursions
(see e.g. \cite{LGM}, or Section \ref{sna-tra} below for the definition of this coding). We write 
$\rho$ for the root of $\t_\zeta$, and we note that $\t_\zeta$ is canonically equipped with 
a volume measure $\mathrm{vol}(\cdot)$. We then consider
Brownian motion indexed by $\t_\zeta$, which we denote by $(V_a)_{a\in\t_\zeta}$ --- we sometimes also
say that $V_a$ is a Brownian label assigned to $a$. Informally,
conditionally on $\t_\zeta$,
$(V_a)_{a\in\t_\zeta}$ is just the centered Gaussian process such that $V_\rho=0$, and 
$\E[(V_a-V_b)^2]=d_\zeta(a,b)$ for every $a,b\in\t_\zeta$, where $d_\zeta$
is the distance on $\t_\zeta$. A formal definition leads to certain technical difficulties 
because the indexing set is random, but these difficulties can be overcome easily using the
formalism of snake trajectories as recalled in Section \ref{sna-tra}. Within this formalism,
the Brownian tree $\t_\zeta$, and the Brownian motion $(V_a)_{a\in\t_\zeta}$ are defined
under a $\sigma$-finite measure $\N_0$ -- see Section \ref{sna-mea} for more details. We note that both the CRT
and Brownian motion indexed by the Brownian tree are important probabilistic objects that appear
as scaling limits for a number of models of combinatorics, interacting particle systems 
and statistical physics (see the introduction of \cite{ALG} for a few related references).
Furthermore, Brownian motion indexed by the Brownian tree is very closely related to
the measure-valued process called super-Brownian motion (see in particular \cite{Zurich}). 

Let us now discuss growth-fragmentation processes, referring to \cite{Ber2} and \cite{BBCK} for additional details. 
The basic ingredient in the construction of a (self-similar) growth-fragmentation process
is a self-similar Markov process  $(X_t)_{t\geq 0}$ with values in $[0,\infty)$ and only negative jumps, which is stopped upon hitting $0$. 
Suppose that $X_0=z>0$, and view $(X_t)_{t\geq 0}$
as the evolution in time of the mass of an initial particle called the Eve particle. 
At each time $t$ where the process $X$ has a jump, we consider that a new particle with mass $-\Delta X_t$ (a
child of the Eve particle)
is born, and the mass of this new particle evolves (from time $t$) again according to the
law of the process $X$, independently of the evolution of the mass of the Eve particle. Then 
each child of the Eve particle has children at discontinuity times of its mass process, and so on.
We consider that a particle dies when its mass vanishes.
Under suitable assumptions (see \cite{Ber2}), we can make sense of the process $(\mathbf{X}(t))_{t\geq 0}$ 
giving for every time $t$ the sequence  (in
nonincreasing order) of masses of all particles alive at that time (if there are only finitely many such particles, the sequence is completed by
adding terms all equal to $0$). The process  $\mathbf{X}$ 
is Markovian and is called the growth-fragmentation
process with Eve particle process $X$. In the preceding description, the process starts from $(z,0,0,\ldots)$, but we can 
get a more general initial value by considering infinitely many Eve particles that evolve independently --- some assumption
is needed on the initial values of these Eve particles so that at every time $t$ the masses of the particles alive can be ranked
in a nonincreasing sequence.

\newpage
\begin{theorem}
\label{main1}
Almost everywhere under the measure $\N_0$, for every $r\geq 0$ and for every connected component $\cc$
of the open set $\{a\in\t_\zeta: V_a>r\}$, the limit 
$$|\partial \cc|:=\lim_{\ve \to 0} \ve^{-2}\mathrm{vol}(\{a\in\cc : V_a<r+\ve\})$$
exists in $(0,\infty)$ and is called the boundary size of $\cc$. For every $r\geq 0$, let $\mathbf{X}(r)$ denote the 
sequence of boundary sizes of all connected components of $\{a\in\t_\zeta: V_a>r\}$ ranked in nonincreasing order. Then, under $\N_0$,
the process $(\mathbf{X}(r))_{r\geq 0}$ is a growth-fragmentation process whose 
Eve particle process $(X_t)_{t\geq 0}$ can be described as follows. 
The process $(X_t)_{t\geq 0}$ is the
self-similar Markov process  with index $\frac{1}{2}$ which in the case 
$X_0=1$ can be represented as
$$X_t=\exp(\xi(\chi(t))),$$
where $(\xi(s))_{s\geq 0}$ is the L\'evy process with only negative jumps and Laplace exponent
\begin{equation}
\label{formula-psi}
\psi(\lambda) = \sqrt{\frac{3}{2\pi}}\,\Big(-\frac{8}{3}\,\lambda + \int_{-\log 2}^0 (e^{\lambda y}-1-\lambda(e^y-1))\,e^{-3y/2}\,(1-e^y)^{-5/2}\,\dd y\Big),
\end{equation}
and $(\chi(t))_{t\geq 0}$ is the time change
\begin{equation}
\label{formula-chi}
\chi(t)=\inf\Big\{s\geq 0: \int_0^s e^{\xi(v)/2}\,\dd v >t\Big\}.
\end{equation}
\end{theorem}

In the setting of Theorem \ref{main1}, we consider the infinite measure $\N_0$, but the statement
still makes sense by conditioning on the initial value $\mathbf{X}(0)$.
The representation of the self-similar Markov process $X$ in terms 
of the L\'evy process $\xi$ is the classical Lamperti representation
of self-similar Markov processes \cite{Lam}. We note that the process $\xi$ drifs to $-\infty$ and $\chi(t)=\infty$ for $t\geq H_0:=\int_0^\infty e^{\xi(v)/2}\,\dd v$,
which simply means that $X_t$ is absorbed at $0$ at time $H_0$. 

\smallskip
It is interesting to relate the growth-fragmentation process of Theorem \ref{main1} to the local times 
of the process $(V_a)_{a\in\t_\zeta}$.
It is known \cite{BMJ} (see also \cite{Sug} for closely related results concerning super-Brownian motion)  that there exists, $\N_0(\dd\omega)$ a.e., a continuous function $(\mathcal{L}_x,x\in\R)$
such that, for every nonnegative measurable function $f$ on $\R$,
$$\int_{\t_\zeta} \mathrm{vol}(\dd a)\,f(V_a) = \int_{\R} \dd x\,f(x)\,\mathcal{L}_x,$$
and we call $\mathcal{L}_x$ the local time at level $x$. Then, for every $r>0$, if $\mathcal{N}_\ve^r$ denotes the 
number of connected components
of $\{a\in\t_\zeta: V_a>r\}$ with boundary size greater than $\ve$, Proposition \ref{local-approx} below gives
$$\ve^{3/2}\,\mathcal{N}^r_\ve\build{\la}_{\ve\to 0}^{}\frac{1}{\sqrt{6\pi}}\,\mathcal{L}_r\;,\quad \N_0\hbox{ a.e.}$$
In other words the suitably rescaled number of fragments of $\mathbf{X}(r)$ with size greater than $\ve$ converges to $\mathcal{L}_r$. 
As a side remark, one might expect the process $(\mathcal{L}_r)_{r\geq 0}$ to be Markovian under $\N_0$, by analogy with the classical Ray-Knight theorems
for local times of linear Brownian motion. This is not the case, but the previous display shows that $\mathcal{L}_r$
is a function of the Markov process $\mathbf{X}(r)$ which obviously contains more information than the local time.

\smallskip
Thanks to the excursion theory developed in \cite{ALG}, we can in fact deduce Theorem \ref{main1} from a simpler statement
valid under the ``positive Brownian snake excursion measure'' $\N^{*}_0$ introduced and studied in \cite{ALG}. 
We refer to Section \ref{sna-mea} for more details, but note that
we can still make sense of the ``genealogical tree'' $\t_\zeta$ and the ``labels'' $V_a$, $a\in\t_\zeta$ under $\N^*_0$. However, we now have 
$V_a\geq 0$ for every $a\in\t_\zeta$, and more precisely the labels $V_b$ are positive along the ancestral line of $a$, except
at the root and possibly at $a$. Informally the measure $\N^*_0$ describes the subtree and the labels corresponding 
under $\N_0$ to a connected component of the set of points with positive labels. One can make sense under $\N^*_0$
of the boundary size
$$\z^*_0:=\lim_{\ve \to 0} \ve^{-2}\mathrm{vol}(\{a\in\t_\zeta : V_a<\ve\})\;,\qquad \N^*_0\hbox{ a.e.}$$
and define the conditional probability measures $\N^{*,z}_0(\cdot)= \N^{*}_0(\cdot\mid \z^*_0=z)$ for every $z>0$.

\begin{theorem}
\label{main2}
Let $z>0$. Almost surely under the measure $\N_0^{*,z}$, for every $r\geq 0$ and for every connected component $\cc$
of the open set $\{a\in\t_\zeta: V_a>r\}$, the limit 
$$\lim_{\ve \to 0} \ve^{-2}\mathrm{vol}(\{a\in\cc : V_a<r+\ve\})$$
exists in $(0,\infty)$ and is called the boundary size of $\cc$. For every $r\geq 0$, let $\mathbf{Y}(r)$ denote the 
sequence of boundary sizes of all connected components of $\{a\in\t_\zeta: V_a>r\}$ ranked in nonincreasing order. Then, under $\N^{*,z}_0$,
the process $(\mathbf{Y}(r))_{r\geq 0}$ is distributed as the growth-fragmentation process of Theorem \ref{main1}
with initial value $\mathbf{Y}(0)=(z,0,0,\ldots)$.
\end{theorem}

Theorem \ref{main1} is a straightforward consequence of Theorem \ref{main2} and the excursion theory of \cite{ALG}. Let us
explain this. Under $\N_0$, the connected components $\cc_1,\cc_2,\ldots$ of $\{a\in\t_\zeta:V_a>0\}$, and the labels
on these components can be represented by snake trajectories $\omega_1,\omega_2,\ldots$ (see Section \ref{sna-tra} for the
definition of snake trajectories). By \cite[Theorem 4]{ALG}, conditionally on the boundary sizes $(|\partial\cc_1|,|\partial\cc_2|,\ldots)$, 
$\omega_1,\omega_2,\ldots$ are independent and the conditional distribution of $\omega_i$
is $\N^{*,|\partial\cc_i]}_0$. In the notation of Theorem \ref{main1}, $\mathbf{X}(0)$ is just 
the (ranked) sequence $(|\partial\cc_1|,|\partial\cc_2|,\ldots)$, and we get that, conditionally on $\mathbf{X}(0)$,
the process $(\mathbf{X}(r))_{r\geq 0}$ is obtained by superimposing {\it independent} processes $(\mathbf{Y}_i(r))_{r\geq 0}$
such that, for every $i\geq 1$,  $\mathbf{Y}_i$ is a growth-fragmentation process started from $(|\partial\cc_i|,0,0,\ldots)$ (by Theorem \ref{main2}).

Theorems \ref{main1} and \ref{main2} have direct applications to the models of random geometry known as the
Brownian map and the Brownian disk. Recall that the Brownian map is a random compact metric space homeomorphic to
the sphere $\mathbb{S}^2$, which is the scaling limit of various classes of random planar maps 
equipped with the graph distance (see in particular \cite{Uniqueness, Mie}). Similarly, the Brownian disk is a random compact metric space homeomorphic to the
closed unit disk of the plane, which appears as the scaling limit of rescaled Boltzmann quadrangulations with a boundary, when the 
size of the boundary grows to infinity (see \cite{Bet,BM,GM2}). We note 
that the papers \cite{Bet,BM} consider Brownian disks with fixed boundary size and volume, but in
the present work we will be interested in the free Brownian disk \cite[Section 1.5]{BM} which has a fixed boundary size but a random volume.
Let us write $\D_z$ for the free Brownian disk with boundary size $z>0$. The space $\D_z$ is equipped with a volume measure denoted by $\mathbf{V}(\dd x)$. The boundary 
$\partial\D_z$ may be defined as the set of all points of $\D_z$
that have no neighborhood homeomorphic to the open unit disk, and for every $x\in\D_z$, we write $H(x)$ for the
``height'' of $x$, meaning the distance from $x$ to the boundary $\partial \D_z$.

\begin{theorem}
\label{main3}
Almost surely, for every $r\geq 0$, for every connected component $\cc$ of $\{x\in\D_z: H(x)>r\}$, the limit
$$|\partial \cc|:= \lim_{\ve\to 0} \ve^{-2}\,\mathbf{V}(\{x\in\cc:H(x)<r+\ve\})$$
exists and is called the boundary size of $\cc$. For every $r\geq 0$, let $\mathbf{Z}(r)$ denote 
 the 
sequence of boundary sizes of all connected components of $\{x\in\D_z: H(x)>r\}$ ranked in nonincreasing order. 
Then, 
the process $(\mathbf{Z}(r))_{r\geq 0}$ is distributed as the growth-fragmentation process of Theorem \ref{main1}
with initial value $\mathbf{Z}(0)=(z,0,0,\ldots)$.
\end{theorem}

As the similarity between the two statements suggests, Theorem \ref{main3} is closely related to Theorem \ref{main2}, and
in fact can be derived from the latter result thanks to the construction of the free Brownian disk (with boundary size $z$) from a snake trajectory 
distributed according to $\N^{*,z}_0$, which is developed in \cite{Disks}. Similarly, we could use Theorem \ref{main1} to derive
a result analogous to Theorem \ref{main3} for the free Brownian map, thanks to the construction of the latter metric space
from a snake trajectory distributed according to $\N_0$ (see e.g.~\cite[Section 3]{Disks}). Rather than writing down 
this statement about the free Brownian map, we give in Section \ref{sec:complements}  an analog of Theorem
\ref{main3} for the Brownian plane, which is an infinite-volume version of the Brownian map that has been shown \cite{Bud,Plane} to be the 
universal scaling limit of infinite random lattices such as the UIPT or the UIPQ. Theorem \ref{growth-frag-plane} below shows that the collection of boundary sizes of 
the connected components of the complement of the ball of radius $r$ centered at the root of the Brownian plane evolves like
the same growth-fragmentation process with indefinite growth starting from $0$ (see \cite[Section 4.2]{BBCK} for a thorough discussion
of this process). 

We next state another result for the Brownian disk, which is 
closely related to Theorem \ref{main3}.

\begin{theorem}
\label{main4}
Let $r>0$. On the event $\{\sup\{H(x):x\in\D_z\}>r\}$, let $\cc_1,\cc_2,\cc_3,\ldots$ be the 
connected components of $\{x\in\D_z: H(x)>r\}$ ranked in nonincreasing order of their boundary sizes,
and for every $j=1,2,\ldots$, let $d_j$ denote the intrinsic metric induced by the Brownian disk metric
on the open set $\D_j$. Then, a.s. on the event $\{\sup\{H(x):x\in\D_z\}>r\}$, for every $j=1,2,\ldots$
the metric $d_j$ has a continuous extension to the closure $\ov{\cc}_j$ of $\cc_j$ in $\D_z$,
and this extension is a metric on $\ov{\cc}_j$. Furthermore, conditionally on the sequence of boundary
sizes $(|\partial\cc_1|,|\partial\cc_2|,\ldots)$, the metric spaces $(\ov{\cc}_1,d_1),(\ov{\cc}_2,d_2),\ldots$
are independent free Brownian disks with respective boundary sizes $|\partial\cc_1|,|\partial\cc_2|,\ldots$.
\end{theorem}

In the same way as Theorem \ref{main3} follows from Theorem \ref{main2}, Theorem \ref{main4} is a 
consequence of a statement (Theorem \ref{law-excu-abo} below) that describes the conditional distribution of the snake trajectories corresponding to
the ``excursions above level $h$'' under $\N^{*,z}_0$, conditionally on the boundary sizes of these excursions.
One might expect that Theorem \ref{law-excu-abo}, which essentially corresponds to the branching property of
growth-fragmentation processes, would be a basic tool for the proof  of Theorem \ref{main2}, but in fact
our proof of Theorem \ref{main2} does not use this branching property. We also note that Theorem \ref{main4} is a Brownian disk 
analog of a result of \cite{Disks} showing that the connected components of the complement of a ball
in the Brownian map are independent Brownian disks conditionally on their volumes and boundary sizes. 

Let us finally mention an interesting corollary of our results.

\begin{corollary}
\label{asymp-extinction}
There exist positive constants $\mathbf{c}_1$ and $\mathbf{c}_2$ such that, for every $r\geq 1$,
$$\mathbf{c}_1\,r^{-6}\leq \N^{*,1}_0\Big(\sup_{a\in\t_\zeta} V_a > r\Big) = \P\Big(\sup_{x\in\D_1} H(x)> r\Big) \leq \mathbf{c}_2\,r^{-6}.$$
\end{corollary}

Corollary \ref{asymp-extinction} immediately follows from Theorem \ref{main2} and Theorem \ref{main3}, by
using the asymptotics for the extinction time of growth-fragmentation processes found in \cite[Corollary 4.5]{BBCK}. 

The proof of Theorem \ref{main2} occupies much of the remaining part of the paper. Let us briefly outline the main steps of this proof. For every $a\in\t_\zeta$ such that $V_a>0$, one can define a 
function $(Z^{(a)}_r)_{0\leq r<V_a}$ such that, for every $r\in[0,V_a)$, $Z^{(a)}_r$ is the boundary size of the
connected component of $\{b\in\t_\zeta: V_b>r\}$ that contains $a$ (see Proposition \ref{pro-exit} below). The function $r\mapsto Z^{(a)}_r$ 
is c\`adl\`ag (right-continuous with left limits) with only negative jumps, and every discontinuity time $r_0$ of this function corresponds to
a ``splitting'' of the connected component containing $a$ into two components, namely the one containing $a$, which
has boundary size $Z^{(a)}_{r_0}$,
and another one with boundary size $|\Delta Z^{(a)}_{r_0}|$. It turns out (Proposition \ref{loc-largest}) that there exists
a unique $\ab\in\t_\zeta$, called the terminal point of the locally largest evolution,  such that, for every
discontinuity time $r_0$ of $r\mapsto Z^{(\ab)}_r$, we have $Z^{(\ab)}_{r_0}> |\Delta Z^{(\ab)}_{r_0}|$ (meaning that
$\ab$ ``stays'' in the component with the larger boundary size) and $V_{\ab}$ is maximal among the labels 
of points satisfying the latter property. Furthermore, the distribution of $(Z^{(\ab)}_r)_{0\leq r<V_\ab}$
is the law of the process $X$ of Theorem \ref{main1} up to its hitting time of $0$ (Proposition \ref{law-loc-largest}).
The process $(Z^{(\ab)}_r)_{0\leq r<V_\ab}$ thus plays the role of the evolution of the mass of the Eve particle.
Furthermore, one verifies that, conditionally on $(Z^{(\ab)}_r)_{0\leq r<V_\ab}$, for every discontinuity time $r_0$, the
connected component that splits off the one containing $\ab$ at time $r_0$ is represented by a
snake trajectory distributed according
to $\N^{*,|\Delta Z^{(a)}_{r_0}|}_0$ (Proposition \ref{excu-loc-largest}).  This provides the recursive structure needed to identify the process
$\mathbf{Y}(r)$ of Theorem \ref{main2} as a growth-fragmentation process.

We finally mention a few recent papers that are related to the present work. We refer
to \cite{Ber2,BBCK,Shi} for the theory of growth-fragmentation processes. As we already mentioned,
Theorem \ref{main3} can be viewed as a continuous version of the main result of \cite{BCK}.
In addition to \cite{Bet,BM,GM2}, free
Brownian disks are discussed in the paper  \cite{MS}, which develops an axiomatic characterization of the Brownian map
as part of a program aiming to equip with
the Brownian map with a canonical conformal structure. 
Brownian disks also play an important role in the recent papers
\cite{GM0,GM3} of Gwynne and Miller motivated by the study of 
statistical physics models on random planar maps. 
Finally we observe that there is an interesting analogy between Theorem \ref{main1} and the fragmentation process
occurring when cutting the CRT at a fixed height. According to \cite{Ber0}, the sequence
of volumes 
of the connected components of the complement of the ball of radius $r$ centered at the root in the CRT
is a self-similar fragmentation process whose dislocation measure has
the form $(2\pi)^{-1/2} (x(1-x))^{-3/2}\,\dd x$. Notice that the L\'evy measure of the process $\xi$ of Theorem \ref{main1}
is the push forward of the measure $\mathbf{1}_{[1/2,1]}(x)\,\sqrt{3/2\pi}\,(x(1-x))^{-5/2}\,\dd x$ under
the mapping $x\mapsto \log x$. 

The paper is organized as follows. Section \ref{sec:preli} gives a number of preliminaries. In particular, we recall the
formalism of snake trajectories, which provides a convenient set-up for the study of Brownian motion indexed by
the Brownian tree. We also give a ``re-rooting'' representation of the measure $\N^{*,z}_0$, which is a key tool in several
subsequent proofs. Section \ref{sec:cc} discusses the connected components of the tree $\t_\zeta$ above a fixed level
and also the components ``above the minimum'': the independence and distributional properties of the latter 
have been studied already in the paper \cite{ALG} and play a basic role in the proof of Theorem \ref{main2}.
Section \ref{sec:bdry-size} is devoted to the existence and properties of the boundary size processes
$(Z^{(a)}_r)_{0\leq r<V_a}$. In this section, we rely on the theory of exit measures for the
Brownian snake \cite{Zurich}. Section \ref{sec:local-larg} introduces the locally largest evolution, and Section \ref{sec:law-loc}
identifies the law of the associated boundary size process (Proposition \ref{law-loc-largest}). A key tool for this identification is Proposition \ref{law-exit-p},
which gives the distribution under $\N_0$ of the exit measure process time-reversed at its last visit to $z>0$.
Section \ref{sec:excu-lle} studies the excursions from the locally largest evolution. Roughly speaking, this study 
provides the recursive structure that shows that the ``children'' of the Eve particle evolve according to
the same Markov process. Theorem \ref{main2} is then proved in Section \ref{gro-frag}, and Theorem \ref{main3}
is derived from Theorem \ref{main2} in Section \ref{sec:slicing}. Section \ref{sec:law-compo} gives the proof of Theorem \ref{main4}.
Finally, Section \ref{sec:complements} contains some complements. In particular, we provide a direct
derivation of the cumulant function associated with our growth-fragmentation processes, which is independent
of the proof of the main results. We also discuss the analog of Theorem \ref{main3} for the Brownian plane, and we 
investigate the relations between local times of $(V_a)_{a\in\t_\zeta}$ and the growth-fragmentation process of Theorem \ref{main1}.
The Appendix gives the proof of two technical results.

\section{Preliminaries}
\label{sec:preli}

\subsection{Snake trajectories}
\label{sna-tra}

Most of this work is devoted to the study of random processes indexed by continuous random trees.
The formalism of snake trajectories, which has been introduced in \cite{ALG}, provides a convenient
framework for this study, and we recall the main definitions that will be needed below.

A (one-dimensional) finite path $\w$ is a continuous mapping $\w:[0,\zeta]\la\R$, where the
number $\zeta=\zeta_{(\w)}$ is called the lifetime of $\w$. We let 
$\W$ denote the space of all finite paths, which is a Polish space when equipped with the
distance
$$d_\W(\w,\w')=|\zeta_{(\w)}-\zeta_{(\w')}|+\sup_{t\geq 0}|\w(t\wedge
\zeta_{(\w)})-\w'(t\wedge\zeta_{(\w')})|.$$
The endpoint or tip of the path $\w$ is denoted by $\wh \w=\w(\zeta_{(\w)})$.
We set $\W_0=\{\w\in\W:\w(0)=0\}$. The trivial element of $\W_0$ 
with zero lifetime is identified with the point $0$ of $\R$. 
Occasionally we will use the notation $\underline\w=\min\{\w(t):0\leq t\leq \zeta_{(\w)}\}$.

\begin{definition}
\label{def:snakepaths}
A snake trajectory (with initial point $0$) is a continuous mapping $s\mapsto \omega_s$
from $\R_+$ into $\W_0$ 
which satisfies the following two properties:
\begin{enumerate}
\item[\rm(i)] We have $\omega_0=0$ and the number $\sigma(\omega):=\sup\{s\geq 0: \omega_s\not =0\}$,
called the duration of the snake trajectory $\omega$,
is finite (by convention $\sigma(\omega)=0$ if $\omega_s=0$ for every $s\geq 0$). 
\item[\rm(ii)] For every $0\leq s\leq s'$, we have
$\omega_s(t)=\omega_{s'}(t)$ for every $t\in[0,\displaystyle{\min_{s\leq r\leq s'}} \zeta_{(\omega_r)}]$.
\end{enumerate} 
\end{definition}

We will write $\S_0$ for the set of all snake trajectories. If $\omega\in \S_0$, we often write $W_s(\omega)=\omega_s$ and $\zeta_s(\omega)=\zeta_{(\omega_s)}$
for every $s\geq 0$. The set $\S_0$ is equipped with the distance
$$d_{\S_0}(\omega,\omega')= |\sigma(\omega)-\sigma(\omega')|+ \sup_{s\geq 0} \,d_\W(W_s(\omega),W_{s}(\omega')).$$
A snake trajectory $\omega$ is completely determined 
by the knowledge of the lifetime function $s\mapsto \zeta_s(\omega)$ and of the tip function $s\mapsto \wh W_s(\omega)$: See \cite[Proposition 8]{ALG}.

Let $\omega\in \S_0$ be a snake trajectory and $\sigma=\sigma(\omega)$. The lifetime function $s\mapsto \zeta_s(\omega)$ codes a
compact $\R$-tree, which will be denoted 
by $\t_\zeta$ and called the {\it genealogical tree} of the snake trajectory. This $\R$-tree is the quotient space $\t_\zeta:=[0,\sigma]/\!\sim$ 
of the interval $[0,\sigma]$
for the equivalence relation
$$s\sim s'\ \hbox{if and only if }\ \zeta_s=\zeta_{s'}= \min_{s\wedge s'\leq r\leq s\vee s'} \zeta_r,$$
and $\t_\zeta$ is equipped with the distance induced by
$$d_\zeta(s,s')= \zeta_s+\zeta_{s'}-2 \min_{s\wedge s'\leq r\leq s\vee s'} \zeta_r.$$
(notice that $d_\zeta(s,s')=0$ if and only if $s\sim s'$, and see e.g.~\cite[Section 3]{LGM} for more information about the
coding of $\R$-trees by continuous functions).  Let $p_\zeta:[0,\sigma]\la \t_\zeta$ stand
for the canonical projection. By convention, $\t_\zeta$ is rooted at the point
$\rho:=p_\zeta(0)=p_\zeta(\sigma)$, and the volume measure $\mathrm{vol}(\cdot)$ on $\t_\zeta$ is defined as the push forward of
Lebesgue measure on $[0,\sigma]$ under $p_\zeta$. For every $a,b\in\t_\zeta$, $\llbracket a,b\rrbracket$ denotes the line segment from
$a$ to $b$, and
the ancestral line of $a$ is the segment $\llbracket \rho,a\rrbracket$ (a point $b$ of $\llbracket \rho,a\rrbracket$ is called an ancestor of $a$, and
we also say that $a$ is a descendant of $b$).
We use the notation $\rrbracket a,b\llbracket$ or $\rrbracket a,b\rrbracket$ with an obvious meaning. Branching points of
$\t_\zeta$ are points $c$ such that $\t_\zeta\backslash\{c\}$ has at least $3$ connected components.

Let us now make a crucial observation: By property (ii) in the definition of  a snake trajectory, the condition $p_\zeta(s)=p_\zeta(s')$ implies that 
$W_s(\omega)=W_{s'}(\omega)$. So the mapping $s\mapsto W_s(\omega)$ can be viewed as defined on the quotient space $\t_\zeta$
(this is indeed the main motivation for introducing snake trajectories: replacing mappings defined on trees, which later will be random trees, by mappings
defined on intervals of the real line).
For $a\in \t_\zeta$, we set $V_a(\omega):=\wh W_s(\omega)$ whenever  $s\in[0,\sigma]$ is such that $a=p_\zeta(s)$ --- by the previous observation this does not
depend on the choice of $s$. We interpret $V_a$ as a ``label'' assigned to the ``vertex'' $a$ of $\t_\zeta$. 
Notice that the mapping $a\mapsto V_a$ is continuous on $\t_\zeta$. 

We will use the notation 
\begin{align*}
&W_*:=\min\{W_s(t): s\geq 0,t\in[0,\zeta_s]\}= \min\{V_a:a\in\t_\zeta\},\\
&W^*:=\max\{W_s(t): s\geq 0,t\in[0,\zeta_s]\}= \max\{V_a:a\in\t_\zeta\}.
\end{align*}

Finally, we will use the notion of a subtrajectory. Let $\omega\in\S_0$ and assume that
the mapping $s\mapsto \zeta_s(\omega)$ is not constant on any nontrivial subinterval
of $[0,\sigma]$ (this will always hold in our applications). Let $a\in\t_\zeta\backslash\{\rho\}$ such that $a$ has strict descendants
and $a$ is not a branching point. Then there exist two times $s_1<s_2$ in $(0,\sigma)$ such that $p_\zeta(s_1)=p_\zeta(s_2)=a$, and the
set $p_\zeta([s_1,s_2])$ consists of all descendants of $a$ in $\t_\zeta$. We define a
new snake trajectory $\omega'$ with duration $s_2-s_1$ by setting, for every $s\geq 0$,
$$\omega'_s(t):=\omega_{(s_1+s)\wedge s_2}(\zeta_{s_1}+t)-\wh\omega_{s_1}\,,\quad \hbox{for } 0\leq t\leq \zeta'_s:=\zeta_{(s_1+s)\wedge s_2}- \zeta_{s_1}.$$
We call $\omega'$ the subtrajectory of $\omega$ rooted at $a$. Informally, $\omega'$ represents the subtree of descendants of $a$ and the associated labels.

\subsection{Re-rooting and truncation of snake trajectories}

We now introduce two important operations on snake trajectories in $\S_0$. The first one 
is the re-rooting operation on $\S_0$ (see \cite[Section 2.2]{ALG}). Let $\omega\in \S_0$ and
$r\in[0,\sigma(\omega)]$. Then $\omega^{[r]}$ is the snake trajectory in $\S_0$ such that
$\sigma(\omega^{[r]})=\sigma(\omega)$ and for every $s\in [0,\sigma(\omega)]$,
\begin{align*}
\zeta_s(\omega^{[r]})&= d_\zeta(r,r\oplus s),\\
\wh W_s(\omega^{[r]})&= \wh W_{r\oplus s}-\wh W_r,
\end{align*}
where we use the notation $r\oplus s=r+s$ if $r+s\leq \sigma$, and $r\oplus s=r+s-\sigma$ otherwise. 
By a remark following the definition of snake trajectories, these prescriptions completely determine $\omega^{[r]}$.

We will write $\zeta_s^{[r]}(\omega)=\zeta_s(\omega^{[r]})$ and $W^{[r]}_s(\omega)=W_s(\omega^{[r]})$. 
The tree $\t_{\zeta^{[r]}}$ is then interpreted as the tree $\t_\zeta$ re-rooted at the vertex $p_\zeta(r)$: More precisely,
the mapping $s\mapsto r\oplus s$ induces an isometry from $\t_{\zeta^{[r]}}$
onto $\t_\zeta$, which maps the root of $\t_{\zeta^{[r]}}$ to $p_\zeta(r)$. Furthermore, the vertices
of $\t_{\zeta^{[r]}}$ receive the ``same'' labels as in $\t_\zeta$,
shifted so that the label of the root is still $0$. 

The second operation is the truncation of snake trajectories. For any $\w\in\W_0$ and $y\in\R$, we set
$$\tau_y(\w):=\inf\{t\in[0,\zeta_{(\w)}]: \w(t)=y\}\,,\; \tau^*_y(\w):=\inf\{t\in(0,\zeta_{(\w)}]: \w(t)=y\}$$
with the usual convention $\inf\varnothing =\infty$ (this convention will be in force throughout this work
unless otherwise indicated). Notice that $\tau_y(\w)=\tau^*_y(\w)$ except possibly if $y=0$. 

Let $\omega\in \S_0$ and $y\in \R$. We set, for every $s\geq 0$,
$$\eta_s(\omega)=\inf\Big\{t\geq 0:\int_0^t \mathrm{d}u\,\mathbf{1}_{\{\zeta_{(\omega_u)}\leq\tau^*_y(\omega_u)\}}>s\Big\}$$
(note that the condition $\zeta_{(\omega_u)}\leq\tau^*_y(\omega_u)$ holds if and only if $\tau^*_y(\omega_u)=\infty$ or $\tau^*_y(\omega_u)=\zeta_{(\omega_u)}$).
Then, setting $\omega'_s=\omega_{\eta_s(\omega)}$ for every $s\geq 0$ defines an element $\omega'$ of $\S_0$,
which will be denoted by  ${\rm tr}_y(\omega)$ and called the truncation of $\omega$ at $y$
(see \cite[Proposition 10]{ALG}). The effect of the time 
change $\eta_s(\omega)$ is to ``eliminate'' those paths $\omega_s$ that hit $y$ (at a positive time when $y=0$) and then survive for a positive
amount of time. The genealogical tree of 
${\rm tr}_y(\omega)$ is canonically and isometrically identified 
with the closed subset of $\t_\zeta$ consisting of all $a$ such that
$V_b(\omega)\not =y$ for every strict ancestor $b$ of $a$ (excluding the root when $y=0$). By abuse of notation, we often write ${\rm tr}_y(W)$ instead of ${\rm tr}_y(\omega)$.

\subsection{Measures on snake trajectories}
\label{sna-mea}

We will be interested in two important measures on $\S_0$. First the Brownian snake excursion 
measure $\N_0$ is the $\sigma$-finite measure on $\S_0$ that satisfies the following two properties: Under $\N_0$,
\begin{enumerate}
\item[(i)] the distribution of the lifetime function $(\zeta_s)_{s\geq 0}$ is the It\^o 
measure of positive excursions of linear Brownian motion, normalized so that, for every $\ve>0$,
$$\N_0\Big(\sup_{s\geq 0} \zeta_s >\ve\Big)=\frac{1}{2\ve};$$
\item[(ii)] conditionally on $(\zeta_s)_{s\geq 0}$, the tip function $(\wh W_s)_{s\geq 0}$ is
a centered Gaussian process with covariance function 
$$K(s,s')= \min_{s\wedge s'\leq r\leq s\vee s'} \zeta_r.$$
\end{enumerate}
Informally, the lifetime process $(\zeta_s)_{s\geq 0}$ evolves under $\N_0$ like a Brownian excursion,
and conditionally on $(\zeta_s)_{s\geq 0}$, each path $W_s$ is a linear Brownian path started from $0$, which
is ``erased'' from its tip when $\zeta_s$ decreases and is ``extended'' when $\zeta_s$ increases.
The measure $\N_0$ can be interpreted as the excursion measure away from $0$ for the 
Markov process in $\W_0$ called the Brownian snake.
We refer to 
\cite{Zurich} for a detailed study of the Brownian snake. 
For every $r>0$, we have
$$\N_0(W^*>r)=\N_0(W_*<-r)={\displaystyle \frac{3}{2r^2}}$$ (see e.g. \cite[Section VI.1]{Zurich}).

The following scaling property is often useful. For $\lambda>0$, for every 
$\omega\in \S_0$, we define $\theta_\lambda(\omega)\in \S_0$
by $\theta_\lambda(\omega)=\omega'$, with
$$\omega'_s(t):= \sqrt{\lambda}\,\omega_{s/\lambda^2}(t/\lambda)\;,\quad
\hbox{for } s\geq 0,\;0\leq t\leq \zeta'_s:=\lambda\zeta_{s/\lambda^2}.$$
Then $\theta_\lambda(\N_0)= \lambda\, \N_{0}$. 

Under $\N_0$, the paths $W_s$, $0<s<\sigma$, take both positive and negative values, simply because they
behave like one-dimensional Brownian paths started from $0$. We will now introduce another important measure
on $\S_0$, which is supported on snake trajectories taking only nonnegative values. For $\delta\geq 0$, let $\S_0^{(\delta)}$
be the set of all $\omega\in \S_0$ such that 
$\sup_{s\geq 0}(\sup_{t\in[0,\zeta_s(\omega)]}|\omega_s(t)|)> \delta$. 
Also set
$$\S_0^+=\{\omega\in\S_0: \omega_s(t)\geq 0 \hbox{ for every }s\geq 0,t\in[0,\zeta_s(\omega)]\}\cap \S^{(0)}_0.$$
There exists a $\sigma$-finite measure $\N^*_0$ on
$\S_0$, which is supported on  $\S_0^+$, 
and gives finite mass to the sets $\S_0^{(\delta)}$ for every $\delta>0$, such that
$$\N^*_0(G)=\lim_{\ve\to 0}\frac{1}{\ve}\,\N_0(G(\mathrm{tr}_{-\ve}(W))),$$
for every bounded continuous function $G$ on $\S_0$ that vanishes 
on $\S_0\backslash\S_0^{(\delta)}$ for some $\delta>0$ (see \cite[Theorem 23]{ALG}). 
Under $\N^*_0$, each of the paths $W_s$, $0<s<\sigma$, starts from $0$, then
stays positive during some time interval $(0,\alpha)$, and is stopped immediately
when it returns to $0$, if it does return to $0$. 

One can in fact make sense of the ``quantity'' of paths $W_s$ that return to $0$ under $\N^*_0$.
To this end, one proves that the limit
\begin{equation}
\label{approxz*0}
\z^*_0:=\lim_{\ve\to 0} \frac{1}{\ve^2} \int_0^\sigma \dd s\,\mathbf{1}_{\{\wh W_s<\ve\}}
\end{equation}
exists $\N^*_0$ a.e. See \cite[Proposition 30]{ALG} for a slightly weaker result --- the stronger form
stated above follows from the results of \cite[Section 10]{Disks}. Notice that replacing the 
limit by a liminf in \eqref{approxz*0} allows us to make sense of $\z^*_0(\omega)$ for every $\omega\in\S_0^+$.
The following conditional versions of the
measure $\N^*_0$ play a fundamental role in the present work. According to
\cite[Proposition 33]{ALG}, there exists a unique collection $(\N^{*,z}_0)_{z>0}$ of probability measures on $\S_0^+$
such that:
\begin{enumerate}
\item[\rm(i)] We have
$$\N^*_0=
 \sqrt{\frac{3}{2 \pi}} 
\int_0^\infty \mathrm{d}z\,z^{-5/2}\, \N^{*,z}_0.$$
\item[\rm(ii)] For every $z>0$, $\N^{*,z}_0$ is supported on $\{\z^*_0=z\}$.
\item[\rm(iii)] For every $z,z'>0$, $\N^{*,z'}_0=\theta_{z'/z}(\N^{*,z}_0)$.
\end{enumerate}
Informally, $\N^{*,z}_0=\N^*_0(\cdot\mid \z^*_0=z)$.

\subsection{Exit measures}
\label{exit-mea}

Let $r\in\R$, $r\not =0$. In a way similar to the definition of $\z^*_0$ above, one can make sense of a quantity that measures the number of paths 
$W_s$ that hit level $r$ under $\N_0$. Precisely, the limit
\begin{equation}
\label{formu-exit}
\z_r:=\lim_{\ve \to 0} \frac{1}{\ve} \int_0^\sigma \dd s\,\mathbf{1}_{\{\tau_r(W_s)\leq \zeta_s<\tau_r(W_s)+\ve\}}
\end{equation}
exists $\N_0$ a.e. Furthermore, $\z_r>0$ if and only if $r\in[W_*,W^*]$, $\N_0$ a.e. This definition of $\z_r$
is a particular case of the theory of exit measures, see \cite[Chapter V]{Zurich}. We note that $\z_r$ is $\N_0$ a.e.
equal to a measurable function of the truncated snake $\mathrm{tr}_r(W)$:  When $r<0$, this can be seen by observing that 
$\z_r$ is the a.e. limit of the quantities $\tilde\z^\ve_r$ introduced in Remark (ii) after Proposition
\ref{approx-exit} below. 

We now recall the special Markov property of the Brownian snake under $\N_0$ (see in particular the appendix
of \cite{subor}).  

\begin{proposition}[Special Markov property]
\label{SMP}
Let $(s_i,s'_i)$, $i\in I$ be the connected components of the
open set $\{s\in[0,\sigma]:\tau_r(W_s)<\zeta_s\}$. For every $i\in I$, set  $a_i:=p_\zeta(s_i)=p_\zeta(s'_i)$ and
let $\omega_i$ be the subtrajectory of $\omega$ rooted at $a_i$.
Then, under the probability measure $\N_0(\cdot\mid r\in[W_*,W^*])$, conditionally on $\mathrm{tr}_r(W)$, the point measure
$\sum_{i\in I}\delta_{\omega_i}$ is Poisson with intensity $\z_r\,\N_0(\cdot)$.
\end{proposition}

Let us now explain 
the relations between exit measures and a certain continuous-state branching process. For $\lambda >0$, we set $\phi(\lambda):=\sqrt{8/3}\,\lambda^{3/2}$ (this notation will be in force
throughout this work). The continuous-state branching process with branching mechanism $\phi$,
or in short the $\phi$-CSBP, is the Feller Markov process $X$ in $\R_+$ whose transition kernels are given
by the following Laplace transform,
\begin{equation}
\label{CSBP-tran}
\E[\exp(-\lambda X_t)\mid X_0=x]= \exp\Big(-x\,\Big(\lambda^{-1/2} + t\sqrt{2/3}\Big)^{-2}\Big),
\end{equation}
for every $x,t\geq 0$ and $\lambda>0$. See e.g. \cite[Chapter II]{Zurich} for basic facts about continuous-state branching processses.

For reasons that will appear later, we now concentrate on the variables $\z_r$ with $r<0$. According to \cite[formula (6)]{CLG}, we have, for every $t>0$,
\begin{equation}
\label{Laplace-exit}
\N_0(1-e^{-\lambda\z_{-t}})= \Big(\lambda^{-1/2} + t\sqrt{2/3}\Big)^{-2}.
\end{equation}
Using both the latter formula and the special Markov property, we get that
the process $(\z_{-r})_{r>0}$ is Markovian under $\N_0$
with the transition kernels of the $\phi$-CSBP, with respect to the filtration $(\g_r)_{r>0}$, where
$\g_r$ denotes the $\sigma$-field generated by $\mathrm{tr}_{-r}(W)$ and the $\N_0$-negligible sets (see \cite[Section 2.5]{ALG}, for more details). Although $\N_0$ is an infinite measure, the 
preceding statement makes sense by considering the process $(\z_{-\delta-r})_{r\geq 0}$
under the probability measure $\N_0(\cdot\mid W_*\leq -\delta)$, for every $\delta>0$. As a 
consequence, the process $(\z_{-r})_{r>0}$ has a c\`adl\`ag modification under $\N_0$,
which we consider from now on.

The distribution of $(\z_{-r})_{r>0}$ under $\N_0$ can be interpreted as 
an excursion measure for the $\phi$-CSBP, in the following sense. Let $\alpha >0$, and let
$$\sum_{i\in I} \delta_{\omega_i}$$
be a Poisson measure with intensity $\alpha\,\N_0$. Set $Y_0=\alpha$ and for every $t>0$,
$$Y_t=\sum_{i\in I} \z_{-t}(\omega_i)$$
(note that this is a finite sum since $\N_0(W_*\leq r)<\infty$ if $r<0$). Then the process
$(Y_t)_{t\geq 0}$ is a $\phi$-CSBP started from $\alpha$. It is enough to verify that $Y$ has the desired one-dimensional marginals and to this end we write, for every $t>0$,
$\E[\exp(-\lambda Y_t)]=\exp(-\alpha\N_0(1-e^{-\lambda\z_{-t}}))$ and we use \eqref{Laplace-exit}.

We note that, for every $z>0$,
\begin{equation}
\label{Zcont0}
\lim_{\ve \to 0}\downarrow \N_0\Big(\sup_{0<t\leq \ve} \z_{-t} \geq z\Big) =0.
\end{equation}
Indeed, assuming that this convergence does not hold, the preceding Poisson representation with $\alpha=z/2$ would imply that the   probability 
of the event $\{\sup\{Y_t:0< t\leq \ve\}\geq z\}$ is bounded below by a positive constant independent of $\ve$, which contradicts the
right-continuity of paths of $Y$ at time $0$. It follows from \eqref{Zcont0} that $\z_{-t}\la 0$ as $t\downarrow0$, $\N_0$ a.e.

Exit measures allow us to state the following formula, which relates
the measures $\N^*_0$ and $\N_0$ via a re-rooting procedure.
Let $G$ be a nonnegative measurable function 
on $\S_0$. Then, 
\begin{equation}
\label{re-root-theo}
\N^*_0\Big(\int_0^\sigma \mathrm{d}r\,G(W^{[r]})\Big)
=2\int_{-\infty}^0 \mathrm{d}b\,\N_0\Big(\z_b\,G(\mathrm{tr}_b(W))\Big).
\end{equation}
See \cite[Theorem 28]{ALG}. 

In view of the subsequent developments, it will be important to have a uniform 
approximation of the exit measure process $(\z_r)_{r<0}$ under $\N_0$. This is the goal
of the next proposition.
For $\w\in\W_0$ and $r\in\R$, we use the notation
$$T_r(\w):=\inf\{t\in[0,\zeta_{(\w)}]: \w(t)<r\}.$$

\begin{proposition}
\label{approx-exit}
For $r<0$ and $\ve>0$, set
$$\z^\ve_r:=\ve^{-2}\int_0^\sigma \dd s\,\mathbf{1}_{\{T_r(W_s)=\infty,\wh W_s<r+\ve\}}.$$
Then, for every $\beta>0$,
$$\sup_{r\in(-\infty,-\beta]} |\z^\ve_r-\z_r| \build{\la}_{\ve \to 0}^{} 0\;,\qquad \N_0\hbox{ a.e.}$$
\end{proposition}

\noindent{\it Remarks.} (i)  The process $(\z_r)_{r<0}$ is c\`agl\`ad (left-continuous with right limits), and the same is true for the process $(\z^\ve_r)_{r<0}$
for every $\ve >0$: If $r_n\uparrow r<0$, we have 
$\mathbf{1}_{\{T_{r_n}(W_s)=\infty\}} \downarrow \mathbf{1}_{\{T_r(W_s)=\infty\}}$ and 
$\mathbf{1}_{\{\wh W_s<r_n+\ve\}}\uparrow\mathbf{1}_{\{\wh W_s<r+\ve\}}$. 
The right limits of the process $(\z^\ve_r)_{r<0}$ are given by
\begin{equation}
\label{right-lim}
\z^\ve_{r+}=\ve^{-2}\int_0^\sigma \dd s\,\mathbf{1}_{\{\un{W}_s>r,\wh W_s\leq r+\ve\}}.
\end{equation}
(ii) The reader may notice that a slightly different approximation is used in \cite[Lemma 14]{ALG}
or in \cite[Proposition 34]{Disks}, where the quantities
$$\wt\z^\ve_r:=\ve^{-2}\int_0^\sigma \dd s\,\mathbf{1}_{\{\zeta_s\leq \tau_r(W_s),\wh W_s<r+\ve\}}.$$
are considered.
If $r$ is fixed, this makes no difference since $\wt\z^\ve_r=\z^\ve_r$ for every $\ve>0$, $\N_0$ a.e. (we may have
$\wt\z^\ve_r\not =\z^\ve_r$ only if $r$ is a local minimum of one of the paths $W_s$, and this occurs with zero $\N_0$-measure).
The point in using $\z^\ve_r$ rather than  $\wt\z^\ve_r$ is the fact that we want a uniform approximation
of $(\z_r)_{r<0}$ and to this end we are looking for c\`agl\`ad approximating processes, which is the case for $r\mapsto \z^\ve_r$ 
but not for $r\mapsto\wt\z^\ve_r$. 

\smallskip
We postpone the proof of Proposition \ref{approx-exit} to the Appendix below. We note that, for every fixed value of $r<0$, the convergence 
$\z^\ve_r \la \z_r$, $\N_0$ a.e., follows from \cite[Proposition 34]{Disks}. Unfortunately, the uniform convergence stated in the proposition
requires more work. 

\subsection{A representation for the measure $\N^{*,z}_0$}
For every $z>0$, set
$$L_z:=\inf\{r<0: \z_r=z\}$$
if $\{r<0: \z_r=z\}$ is not empty, and $L_z=0$ otherwise. 

\begin{lemma}
\label{law-Lz}
We have
\begin{equation}
\label{law-L}
\N_0(L_z<0)= \frac{1}{2z}.
\end{equation}
\end{lemma}

\proof The fact that $\N_0(L_z<0)=C/z$ for some positive constant $C$
is easy by a scaling argument, but we need another argument to get the value of $C$. Let $\ve >0$. We have
$$
\N_0(L_z<0)=\N_0\Big(\sup_{t>0} \z_{-t} \geq z\Big)
=\N_0\Big(\sup_{0<t\leq \ve} \z_{-t} \geq z\Big) +\N_0\Big(\mathbf{1}_{\{\sup_{0<t\leq \ve} \z_{-t}< z\}}\mathbf{P}_{\z_{-\ve}}\Big(\sup_{t\geq 0}\mathscr{Y}_t \geq z\Big)\Big),
$$
where we use the notation $(\mathscr{Y}_t)_{t\geq 0}$ for a $\phi$-CSBP that starts from $x$ under the probability measure $\mathbf{P}_x$, for every $x\geq 0$. 
By the classical Lamperti representation for CSBPs \cite{Lam0,CLU}, $(\mathscr{Y}_t)_{t\geq 0}$ can be written as a time change of a stable L\'evy process
with index $3/2$ and no negative jumps. The explicit solution of the two-sided exit problem for such L\'evy processes (see
\cite[Theorem VII.8]{Ber}) now gives 
$$\mathbf{P}_{\z_{-\ve}}\Big(\sup_{t\geq 0}\mathscr{Y}_t \geq z\Big)= 1-\sqrt{\Big(1-\frac{\z_{-\ve}}{z}\Big)^+}.$$
Using also \eqref{Zcont0}, we get that
 $$\N_0(L_z<0)= \N_0\Bigg( 1-\sqrt{\Big(1-\frac{\z_{-\ve}}{z}\Big)^+}\Bigg) + o(1)\;,\quad\hbox{as }\ve\to 0.$$
 
For every $\delta>0$, $\N_0(\mathbf{1}_{\{\z_{-\ve} >\delta\}}\z_{-\ve})=o(1)$ as $\ve\to 0$. Indeed, this follows from
the formula $\N_0(\z_\ve(1-e^{-\lambda \z_{-\ve}}))=1-(1+\sqrt{2/3}\, \ve\lambda^{1/2})^{-3}$, which is a consequence of \eqref{Laplace-exit}. Thanks to this
observation and to the fact that $\sqrt{1-x}=1-\frac{x}{2}+o(x)$ when $x\to 0$, we find that
$$\N_0(L_z<0)=\frac{1}{2} \N_0\Big(\frac{\z_{-\ve}}{z}\Big) + o(1)$$
as $\ve \to 0$. The lemma follows since $\N_0(\z_{-\ve})=1$ for every $\ve >0$. \endproof

The following proposition, which will play an important role, provides an analog of formula
\eqref{re-root-theo} where $\N^*_0$ is replaced by the conditional measure $\N^{*,z}_0$. 

\begin{proposition}
\label{repre}
For any nonnegative measurable function $G$ on $\S_0$, for every $z>0$,
$$z^{-2}\,\N^{*,z}_0\Big(\int_0^\sigma \dd s\,G(W^{[s]})\Big)= \N_0\Big( G(\mathrm{tr}_{L_z}(W))\,\Big| L_z<0\Big).$$
\end{proposition}

\noindent{\it Remark.} When $G=1$, one recovers the known formula $\N^{*,z}_0(\sigma)=z^2$, see the remark following Proposition 15 in \cite{Disks}. 

\proof We may and will assume that $G$ is bounded and continuous. We use the same notation $(\mathscr{Y}_t,\mathbf{P}_x)$
as in the previous proof and we also set
$\Lambda_z:=\sup\{t\geq 0: \mathscr{Y}_t=z\}$
with the convention $\sup\varnothing=0$.

As a consequence of \eqref{re-root-theo}, the formula
\begin{equation}
\label{re-root-rep-tec0}
\N^*_0\Big(\int_0^\sigma \mathrm{d}r\,\varphi(\z^*_0)\,G(W^{[r]})\Big)
=2\int_{-\infty}^0 \mathrm{d}b\,\N_0\Big(\z_b\,\varphi(\z_b)\,G(\mathrm{tr}_b(W))\Big),
\end{equation}
holds for any nonnegative measurable function 
$\varphi$ on $[0,\infty)$. To derive \eqref{re-root-rep-tec0} from \eqref{re-root-theo}, notice that \eqref{approxz*0} and Proposition \ref{approx-exit}
allow us to write $\z^*_0=\Gamma(W^{[r]})$, $\N^*_0$ a.e., and $\z_b=\Gamma(\mathrm{tr}_b(W))$, $\N_0$ a.e., with the {\it same}
measurable function $\Gamma$ on $\S_0$.

Let us fix $z_0>0$ and a continuous function $\varphi$ on $\R_+$ which is supported on a compact
subset of $(0,\infty)$ and such that $\varphi(z_0)>0$. We observe that, for any $b<0$, we have
\begin{equation}
\label{MarkovZ}
\N_0\Big(\mathbf{1}_{\{b-\ve\leq L_{z_0}<b\}}\,\z_b\,\varphi(\z_b)\,G(\mathrm{tr}_b(W))\Big)
= \N_0\Big(h_\ve(\z_b,z_0)\,\z_b\,\varphi(\z_b)\,G(\mathrm{tr}_b(W))\Big),
\end{equation}
where the function $h_\ve$ is defined for every $z>0$ by
$$h_\ve(z,z_0)=\mathbf{P}_z(0<\Lambda_{z_0}\leq \ve).$$
To get \eqref{MarkovZ}, we use the Markov property  of the process $(\z_{-r})_{r> 0}$ 
(with respect to the filtration $(\g_r)_{r>0}$ introduced 
in Section \ref{exit-mea}) at time $-b$.

By combining \eqref{re-root-rep-tec0} (with $\varphi(z)$ replaced by $h_\ve(z,z_0)\,\varphi(z)$) and \eqref{MarkovZ}, we get
\begin{align}
\label{key00}
\N^*_0\Big(\int_0^\sigma \mathrm{d}r\,h_\ve(\z^*_0,z_0)\,\varphi(\z^*_0)\,G(W^{[r]})\Big)
&=
2\int_{-\infty}^0 \mathrm{d}b\,\N_0\Big(\mathbf{1}_{\{b-\ve\leq L_{z_0}<b\}}\,\z_b\,\varphi(\z_b)\,G(\mathrm{tr}_b(W))\Big)\nonumber\\
&= 2\,\N_0\Big(\int_{L_{z_0}}^{(L_{z_0}+\ve)\wedge 0} \dd b \,\z_b\,\varphi(\z_b)\,G(\mathrm{tr}_b(W))\Big)
\end{align}
Let us multiply the right-hand side of \eqref{key00} by $\ve^{-1}$ and study its limit as $\ve\to 0$. By Lemma 11 in \cite{ALG}
we know that $\mathrm{tr}_b(W)\la \mathrm{tr}_{L_{z_0}}(W)$ as $b\downarrow L_{z_0}$,
$\N_0$ a.e. on $\{L_{z_0}<0\}$. It follows that
\begin{equation}
\label{key02}
\lim_{\ve\to 0} \frac{2}{\ve}\N_0\Big(\int_{L_{z_0}}^{(L_{z_0}+\ve)\wedge 0} \dd b \,\z_b\,\varphi(\z_b)\,G(\mathrm{tr}_b(W))\Big)
= 2z_0\,\varphi(z_0)\,\N_0\Big(G(\mathrm{tr}_{L_{z_0}}(W))\,\mathbf{1}_{\{L_{z_0}<0\}}\Big),\end{equation}
where dominated convergence is easily justified thanks to our assumptions on $\varphi$ and the property $\N_0(L_{z_0}<0)<\infty$.
On the other hand, properties (i) and (ii) stated at the end of Section \ref{sna-mea} allow us to rewrite the left-hand side of \eqref{key00} as
\begin{equation}
\label{key01}
\N^*_0\Big(\int_0^\sigma \mathrm{d}r\,h_\ve(\z^*_0,z_0)\,\varphi(\z^*_0)\,G(W^{[r]})\Big)
= \sqrt{\frac{3}{2\pi}}\,\int_0^\infty \frac{\dd z}{z^{5/2}}\,h_\ve(z,z_0)\,\varphi(z)\,\N^{*,z}_0\Big(\int_0^\sigma \dd r\,G(W^{[r]})\Big).
\end{equation}
Consider the special case $G=1$. We deduce from the convergence \eqref{key02}, using also the
formula $\N^{*,z}_0(\sigma)=z^2$ and the identities \eqref{key00} and
\eqref{key01}, that
\begin{equation}
\label{key03}
\lim_{\ve\to 0} \frac{1}{\ve}\,\sqrt{\frac{3}{2\pi}}\,\int_0^\infty \frac{\dd z}{z^{1/2}}\,h_\ve(z,z_0)\,\varphi(z)= 2z_0\,\varphi(z_0) \N_0(L_{z_0}<0).
\end{equation}
For a general (bounded and continuous) function $G$, a simple
scaling argument shows that the function
$z\mapsto z^{-2}\,\N^{*,z}_0(\int_0^\sigma \dd r\,G(W^{[r]}))$
is also bounded and continuous on $(0,\infty)$. We may thus apply \eqref{key03} with $\varphi(z)$ replaced by the function
$$z\mapsto \varphi(z)\,z^{-2}\N^{*,z}_0\Big(\int_0^\sigma \dd r\,G(W^{[r]})\Big)$$
and we get
\begin{equation}
\label{key04}
\lim_{\ve\to 0} \frac{1}{\ve}\sqrt{\frac{3}{2\pi}}\int_0^\infty \frac{\dd z}{z^{5/2}}h_\ve(z,z_0)\varphi(z)\N^{*,z}_0\Big(\int_0^\sigma \dd r\,G(W^{[r]})\Big)
=\frac{2\varphi(z_0)}{z_0}\,\N^{*,z_0}_0\Big(\int_0^\sigma \dd r\,G(W^{[r]})\Big)\,\N_0(L_{z_0}<\infty)
\end{equation}
From the identities \eqref{key00} and \eqref{key01}, the right-hand sides of \eqref{key02} and \eqref{key04} are equal, which gives
 the desired result. \endproof

\subsection{The exit measure process time-reversed at $L_z$}

 The goal of this section is to prove the following proposition. Recall that for a L\'evy process $\xi$
with only negative jumps we define its Laplace exponent $\psi (\lambda)$ by
$$\E[\exp(\lambda \xi(t))]= \exp(t\psi(\lambda))\;,\quad \lambda\geq 0.$$
We use the notation $\z_{r+}$ for the right limit of $u\mapsto \z_u$ at $r$.

\begin{proposition}
\label{law-exit-p}
Set $\z_r=0$ for $r\geq 0$. Under $\N_0(\cdot\mid L_z<0)$, the process $(\z_{(L_z+r)+})_{r\geq 0}$ is distributed as 
a self-similar Markov process $(X^\circ_r)_{r\geq 0}$ with index $\frac{1}{2}$ starting from $z$, which can be represented as
$$X^\circ_t=z\,\exp(\xi^\circ(\chi^\circ(z^{-1/2}t))),$$
where $(\xi^\circ(s))_{s\geq 0}$ is the L\'evy process with only negative jumps and Laplace exponent
$$\psi^\circ(\lambda)=\sqrt{\frac{3}{2\pi}} \int_{-\infty}^0 (e^{\lambda y}-1-\lambda(e^y-1))\,e^{y/2}\,(1-e^y)^{-5/2}\,\dd y,$$
and $(\chi^\circ(t))_{t\geq 0}$ is the time change
$$\chi^\circ(t)=\inf\Big\{s\geq 0: \int_0^s e^{\xi^\circ(v)/2}\,\dd v >t\Big\}.$$
\end{proposition}

We note that the L\'evy process $\xi^\circ$ drifts to $-\infty$, and the quantity
$$H^\circ_0:=z^{1/2}\int_0^\infty e^{\xi^\circ(v)/2}\,\dd v$$
is finite a.s. For $t\geq H^\circ_0$, we have $\chi^\circ(z^{-1/2}t)=\infty$ and $\xi^\circ(\chi^\circ(z^{-1/2}t))=-\infty$. Thus $H^\circ_0$
is the hitting time of $0$ by $X^\circ$, and $X^\circ$ is absorbed at $0$. 

\proof
Let $(U_t)_{t\geq 0}$ denote a stable L\'evy process with index $3/2$ and no negative jumps, whose distribution is
characterized by the formula
$$\E[\exp(-\lambda U_t)]=\exp(t\phi(\lambda))\;,\quad \lambda>0,\,t\geq 0$$
where $\phi(\lambda)=\sqrt{8/3}\,\lambda^{3/2}$ as previously. If $\un{U}_t:=\min\{U_s:0\leq s\leq t\}$, the
process $U_t-\un{U}_t$ is a strong Markov process for which $0$ is a regular point. Furthermore, $-\un{U}_t$
serves as a local time at $0$ for $U-\un{U}$. We refer to \cite{Ber}, especially Chapters VII and VIII, for these standard facts about L\'evy processes.
We denote the excursion measure of $U-\un{U}$ away from $0$, corresponding to the local time $-\un{U}$, by $\bn$. Then $\bn$
is a $\sigma$-finite measure on the Skorokhod space $\mathbb{D}(\R_+,\R_+)$. 

For notational convenience, we write $\bar \z_x= \z_{-x}$ for $x>0$ and $\bar\z_0=0$. Notice that $\bar\z$ has c\`adl\`ag sample paths. We also set
$$\bar L_z=-L_z=\sup\{x>0: \bar \z_x=z\}$$
with the convention $\sup\varnothing=0$.

\begin{lemma}
\label{Lamperti-excu}
For every $x\geq 0$, set
$$\eta(x):=\inf\{y>0:\int_0^y \bar\z_u\,\dd u >x\}.$$
Let $\y_x=\bar\z_{\eta(x)}$ if $\eta(x)<\infty$ and $\y_x=0$ otherwise.
Then the distribution of $(\y_x)_{x\geq 0}$ under $\N_0$ is $\bn$. 
\end{lemma}

This lemma is basically a version for excursion measures of the Lamperti
representation \cite{Lam0,CLU} connecting continuous-state branching processes with  L\'evy processes.
As we were unable to find a precise reference, we provide a proof in the Appendix below. 

On the event $\{\bar L_z>0\}$, set
$$\Lambda_z:= \sup\{x\geq 0: \y_x=z\}=\int_0^{\bar L_z} \bar \z_s\,\dd s.$$
Still on the event $\{\bar L_z>0\}$, we then introduce the time-reversed processes
$$
\check \z_u=\left\{
\begin{array}{ll}
\bar\z_{(\bar L_z-u)-}\quad&\hbox{if } 0\leq u<\bar L_z,\\
0&\hbox{if }u\geq \bar L_z,
\end{array}\right.$$
and
$$
\check \y_u=\left\{
\begin{array}{ll}
\y_{(\Lambda_z-u)-}\quad&\hbox{if } 0\leq u<\Lambda_z,\\
0&\hbox{if }u\geq \Lambda_z.
\end{array}\right.$$
 We note that we have again the Lamperti representation
\begin{equation}
\label{Lamp0}
\check \z_t = \check \y_{\gamma(t)}\;,\hbox{ with }\gamma(t)=\inf\{u\geq 0: \int_0^u \frac{\dd v}{\check \y_v} >t\}.
\end{equation}

Next as a consequence of Lemma \ref{Lamperti-excu} and Theorem 4 in \cite{Chau}, we know that
the process $(\check\y_u)_{u\geq 0}$ is distributed under $\N_0(\cdot \mid L_z<0)$ as the L\'evy process
$-U$ started from $z$ and conditioned to hit zero continuously before hitting $(-\infty,0)$, and stopped at that hitting time.
We refer to Section 4 of \cite{Chau} for a discussion of the latter process. Furthermore we can then use Corollary 3 of Caballero
and Chaumont \cite{CabCha} to obtain that the process $(\check\y_u)_{u\geq 0}$ under $\N_0(\cdot \mid L_z<0)$
has the distribution of a self-similar Markov process $(X'_u)_{u\geq 0}$ which can be represented in the form
$$X'_u=z\,\exp(\xi^\circ(\chi'(z^{-3/2}u))),$$
where $\xi^\circ$ is the L\'evy process in the statement of the proposition\footnote{In order for the reader to recover the exact form
of the Laplace exponent $\psi^\circ$ in the proposition, we mention the following minor inaccuracy in \cite{CabCha}: In formula
(23)
of the latter paper, $+ c_-$ should be replaced by $-c_-$.}, and $(\chi'(t))_{t\geq 0}$ is the time change
$$\chi'(t)=\inf\Big\{s\geq 0: \int_0^s e^{3\xi^\circ(v)/2}\,\dd v >t\Big\}.$$
(Note that the self-similarity index of $X'$ is $3/2$ as the one for $U$.) 

Recalling \eqref{Lamp0}, we see that $(\check \z_t)_{t\geq 0}$ has the same distribution as $(X'_{\gamma'(t))})_{t\geq 0}$ where
$$\gamma'(t)=\inf\{u\geq 0: \int_0^u \frac{\dd v}{X'_v} >t\}$$
and $X'_\infty=0$ by convention. Let $H'_0:=\inf\{t\geq 0: X'_t=0\}$ and $K'_0:=\int_0^{H'_0} (X'_v)^{-1}\,\dd v$,
so that $\gamma'(t)<H'_0$ if $t<K'_0$ and $\gamma'(t)=\infty$ if $t\geq K'_0$. Simple manipulations show that 
\begin{align*}\chi'(z^{-3/2}\gamma'(t))=\int_0^t \frac{\dd s}{\sqrt{X'_{\gamma'(s)}}}&=\inf\{u\geq 0: z^{1/2}\int_0^u \exp(\frac{1}{2} \xi^\circ(v))\,\dd v >t\}\\
&=\chi^\circ(z^{-1/2}t)
\end{align*}
 if $t<K'_0= z^{1/2} \int_0^\infty e^{\xi^\circ(u)/2}\,\dd u$, whereas 
$\chi'(z^{-3/2}\gamma'(t))=\infty$ if $t\geq K'_0$. In both cases we get
$X'_{\gamma'(t)}= z\,\exp(\xi^\circ(\chi^\circ(z^{-1/2}t)))= X^\circ_t$, with the notation of the proposition. 
We conclude that $(\check \z_t)_{t\geq 0}$ has the same distribution as $(X^\circ_t)_{t\geq 0}$. This is the desired result since by construction $\check \z_t=\z_{L_z+t}$. \hfill$\square$
%
%

%

\section{Special connected components of the genealogical tree}
\label{sec:cc}

\subsection{Components above a level}
\label{comp-level}

In this section and the next one, we formulate certain definitions and facts that make sense for
a deterministic snake trajectory satisfying some regularity properties. 
We fix $\omega\in\S_0$ and consider the associated genealogical tree $\t_\zeta$.
We say that $x\in\R$ is a local minimum of $\omega$ if there exist two distinct points
$a_1,a_2\in \t_\zeta$ and a point $b\in\rrbracket a_1,a_2\llbracket$ such that
$$V_b=\min_{c\in\llbracket a_1,a_2\rrbracket} V_c = x.$$
We then also say that $b$ is a point of local minimum. Clearly the set of all local minima is countable.

We will assume the following regularity properties:
\begin{enumerate}
\item[(i)] local minima are distinct: if $b,b'$ are two distinct points of local minimum, $V_b\not =V_{b'}$;
\item[(ii)] no branching point is a point of local minimum;
\item[(iii)] for every $x\in\R$, $\mathrm{vol}(\{c\in\t_\zeta: V_c=x\})=0$.
\end{enumerate}
All these properties hold $\N_0$ a.e. and $\N^{*,z}_0$ a.e. (for (iii), the case of $\N_0$
follows from the fact that the push forward of $\mathrm{vol}(\dd a)$ under the mapping $a\mapsto V_a$
has a continuous density \cite{BMJ}, and one can then use Proposition \ref{repre}
to deal with $\N^{*,z}_0$). Notice that (i)
implies that the mapping $c\mapsto V_c$ cannot be constant on a nontrivial line segment of $\t_\zeta$.

In the remaining part of this section, we assume in addition that $\omega\in \S^+_0$. We set
$$\t_\zeta^\circ:=\{a\in\t_\zeta:V_a>0\}.$$
Let us  fix $a\in \t_\zeta^\circ$. For every $r\in[0,V_a)$, let $\cc^{(a)}_r$ denote the connected component 
of $\{b\in\t_\zeta: V_b>r\}$ that contains $a$. We note that $\cc^{(a)}_{r'}\subset \cc^{(a)}_{r}$ if $r<r'<V_a$.
Let $\ov{\cc}^{(a)}_r$ stand for the closure of $\cc^{(a)}_r$, and if $r\in(0,V_a)$ set
$$\cc^{(a)}_{r-}:=\bigcap_{r'\in[0,r)} \cc^{(a)}_{r'}.$$
We have always 
$\ov{\cc}^{(a)}_r \subset \cc^{(a)}_{r-}$
and equality holds if and only if $r\notin D^{(a)}$, where the set $D^{(a)}$
is defined by
$$D^{(a)}:=\{r\in (0,V_a): \exists b\in\t_\zeta\backslash\{a\},\;V_b>r\hbox{ and } \min_{c\in\llbracket a,b\rrbracket} V_c =r\}.$$
Note that $D^{(a)}$ is a subset of the set of all local minima. 

If $r\in D^{(a)}$ and $b\not = a$ is such that $V_b>r$ and $\min_{c\in\llbracket a,b\rrbracket} V_c =r$,
then there exists a unique $c_0\in \rrbracket a,b\llbracket$ such that $V_{c_0}=r$, and $c_0$ does not depend on the choice 
of $b$ (because local minima are distinct by (i) above). Note that $c_0$ cannot be a branching point of the tree 
$\t_\zeta$, by (ii). We can then set
$$\check\cc^{(a)}_r= \{b\in\t_\zeta: c_0\in \llbracket a,b\llbracket \hbox{ and } V_c>r\hbox{ for every }c\in\rrbracket c_0,b\rrbracket\},$$
and $\cc^{(a)}_{r-}$ is the closure of the union ${\cc}^{(a)}_r \cup \check\cc^{(a)}_r$. Notice that $\check\cc^{(a)}_r=\cc^{(b)}_r$
for any $b\in \check\cc^{(a)}_r$. For future use, we note that the
boundary of $\cc^{(a)}_r$, or of $\check\cc^{(a)}_r$, has zero volume (by (iii)).

\subsection{Excursions above the minimum}
\label{excu-mini}

Let us consider $\omega\in \S_0$, and assume that the conditions (i)---(iii) of the previous section hold.
Recall our notation $\rho$ for the root of $\t_\zeta$ and note that $V_\rho=0$. In a way very similar to the 
definition of $D^{(a)}$ above we now set
$$D(\omega)=\{r<0: \exists a\in\t_\zeta, V_a>r\hbox{ and  }\min_{c\in\llbracket \rho,a\rrbracket} V_c=r\}.$$
Let us fix $r\in D$. Then $r$ is a local minimum and we let $c_0$ be  the 
uniquely determined point of local minimum such that $V_{c_0}=r$. The same arguments as in the previous section allow us to single out a particular component
of $\{c\in\t_\zeta:V_c>r\}$ by setting
$$\check\cc_r=\{a\in\t_\zeta: c_0\in \llbracket \rho,a\llbracket \hbox{ and } V_c>r\hbox{ for every }c\in\rrbracket c_0,a\rrbracket\}.$$
It is convenient to represent $\check\cc_r$ and the labels on this component by a snake trajectory $\omega^r$, which may be defined
as follows. Since the point $c_0$ has strict descendants in the tree $\t_\zeta$ and is not a branching point, we can make sense 
of the subtrajectory rooted at $c_0$, which we denote by $\tilde \omega^r$ (see Section \ref{sna-tra}). We are in fact only
interested in those descendants of $c_0$ that lie in $\check\cc_r$, and for this reason, we 
consider the truncation $\omega^r=\mathrm{tr}_0(\tilde \omega^r)$. 

Write $\t_\zeta(\omega^r)$ for the genealogical
tree of $\omega^r$ and, as previously, let $\t_\zeta^\circ(\omega^r)$ denote the subset of $\t_\zeta(\omega^r)$ consisting of 
points with positive labels. Then $\check \cc_r$ is identified to $\t_\zeta^\circ(\omega^r)$ via a volume preserving
isometry, in such a way that, for every
$a\in\check\cc_r$, we have $V_a(\omega)=r+V_{\tilde a}(\omega^r)$ if $\tilde a$ is the point of $\t_\zeta^\circ(\omega^r)$ corresponding
to $a$. Consequently, for every $\delta\geq 0$, connected components of $\{a\in\t_\zeta: V_a(\omega)>r+\delta\}$
contained in $\check\cc_r$ are in one-to-one correspondence with connected components of $\{a\in\t_\zeta(\omega^r): V_a(\omega)>\delta\}$. The latter fact 
will be important for our applications in Section \ref{gro-frag} below. 

We call $\omega^r$, $r\in D$, the excursions of 
$\omega$ above the minimum. We refer to \cite[Section 3]{ALG} for a (slightly different) more detailed presentation. 

The following theorem, which is one of the main results of \cite{ALG}, identifies the conditional distribution 
of the excursions $\omega^r$, $r\in D$, under $\N_0$ and conditionally on the exit measure process
$(\z_r)_{r<0}$.

\begin{theorem}{\rm \cite[Proposition 36, Theorem 40]{ALG}}
\label{mainALG}
$\N_0(\dd \omega)$ a.e., $D$ coincides with the set of all discontinuity times of the process $(\z_r)_{r<0}$. 
We can thus write $D=\{r_1,r_2,\ldots\}$ where $r_1,r_2,\ldots$ is the sequence of these
discontinuity times ordered so that $|\Delta\z_{r_1}|>|\Delta\z_{r_2}|>\cdots$.
Then, under $\N_0$ and conditionally on $(\z_r)_{r<0}$, the random snake trajectories $\omega^{r_1},\omega^{r_2},\ldots$
are independent, and for every $i\geq 1$ the conditional distribution of $\omega^{r_i}$
is $\N^{*,|\Delta\z_{r_i}|}_0$.
\end{theorem}

\section{Measuring the boundary size of  components above a level}
\label{sec:bdry-size}

We will now argue under $\N^{*,z}_0$ for a fixed $z>0$. The measure 
$\N^{*,z}_0$ is supported on $\S^+_0$, and so we may use the notation
introduced in Section \ref{comp-level}. Recall in particular that
$\t_\zeta^\circ:=\{a\in\t_\zeta:V_a>0\}$. 

\smallskip
For $a\in\t_\zeta^\circ$, $r\in[0,V_a)$ and $\ve >0$, we set
$$Z^{(a),\ve}_r:= \ve^{-2} \,\vol(\{b\in\cc^{(a)}_r: V_b\leq r+\ve\}).$$

\begin{proposition}
\label{pro-exit}
The following properties hold $\N^{*,z}_0$ a.e.

\noindent{\rm (i)} For every $a\in\t^\circ_\zeta$, $(Z^{(a),\ve}_r)_{r\in[0,V_a)}$ converges as $\ve \to 0$, uniformly on
$[0,V_a-\beta]$ for every $\beta>0$, to a limiting c\`adl\`ag function $(Z^{(a)}_r)_{r\in[0, V_a)}$ with only negative jumps, which takes positive values on $[0,V_a)$ and is such that $Z^{(a)}_0=z$.

\noindent{\rm(ii)} If $a,a'\in \t_\zeta^\circ$, we have $Z^{(a)}_r=Z^{(a')}_r$ for every $r\in [0,\min_{c\in\llbracket a,a'\rrbracket} V_c)$.

\noindent{\rm (iii)} Let $a\in \t_\zeta^\circ$. The set of discontinuities of $r\mapsto Z^{(a)}_r$ is $D^{(a)}$. If $r\in D^{(a)}$
we have
$$Z^{(a)}_{r-}= Z^{(a)}_r + Z^{(b)}_r$$
where $b$ is an arbitrary point of $\check\cc^{(a)}_r$. Moreover $Z^{(a)}_r \not= Z^{(b)}_r$. 
\end{proposition}

\proof Recall from Proposition \ref{approx-exit} the notation $\z^\ve_r(\omega)$ for $\omega\in \S_0$ and $r<0$, and the fact that 
$r\mapsto \z^\ve_r(\omega)$ is c\`agl\`ad. We let $\Theta_z$ be the set of all snake trajectories 
$\omega\in\S_0$ such that $W_*(\omega)<0$ and:
\begin{itemize}
\item[(a)] $(\z^\ve_r(\omega))_{r\in[W_*(\omega),0)}$ converges 
as $\ve\to 0$ to a limiting c\`agl\`ad function  $(\z_r(\omega))_{r\in[W_*(\omega),0)}$, uniformly
on $[W_*(\omega),-\beta]$ for every $\beta\in(0,-W_*(\omega))$;
\item[(b)]
the set of discontinuity times of this limiting
function is $D(\omega)\cap(W_*(\omega),0)$;
\item[(c)] the function $(\z_r(\omega))_{r\in[W_*(\omega),0)}$ takes positive values on $[W_*(\omega),0)$, and takes the value $z$ for
$r=W_*(\omega)$;
\item[(d)]
$\z_r(\omega)\not = |\Delta \z_r(\omega)|$
for every $r\in D(\omega)\cap (W_*(\omega),0)$.
\end{itemize}

It follows from Proposition \ref{approx-exit} and the first assertion of Theorem \ref{mainALG} that $\mathrm{tr}_{L_z}(W)$ belongs to $\Theta_z$,
$\N_0$ a.e. on the event $\{L_z<0\}$. We also use the fact that the process $(\z_{-r}(\omega))_{r>0}$ evolves 
under $\N_0$ as a $\phi$-CSBP (and therefore
as the time change of a stable L\'evy process) to obtain the property $\z_r(\omega)\not = |\Delta \z_r(\omega)|$
when $r\in D(\omega)$. 

Taking $G=\mathbf{1}_{\Theta_z}$ in Proposition \ref{repre}, we also get that, $\N^{*,z}_0(\dd \omega)$ a.s., for $\dd s$ a.e. $s\in[0,\sigma]$,
the re-rooted snake trajectory $W^{[s]}$ belongs to $\Theta_z$. 
So let us fix $\omega\in\S_0$ such that the preceding assertion holds. We can then take a sequence $s_1,s_2,\ldots$ dense 
in $[0,\sigma]$ such that $\omega^{[s_i]}$ belongs to $\Theta_z$ for every $i=1,2,\ldots$. Setting $a_i=p_\zeta(s_i)$, we also know that $a_1,a_2,\ldots$ all belong to $\t_\zeta^\circ$ (otherwise $W_*(\omega^{[s_i]})=0$). We now observe that $W_*(\omega^{[s_i]})=-V_{a_i}(\omega)$, and, for every $r\in[0,V_{a_i}(\omega))$,
$$Z^{(a_i),\ve}_r(\omega)= \z^\ve_{(r-V_{a_i}(\omega))+}(\omega^{[s_i]}).$$
This is a simple consequence of our definitions and formula \eqref{right-lim} for the right limits $\z^\ve_{r+}$. 

Since $\omega^{[s_i]}$ belongs to $\Theta_z$, we deduce from the last display and assertion (a) above that the convergence 
stated in part (i) of the proposition holds when $a=a_i$, and that, for every $r\in[0,V_{a_i}(\omega))$,
$$Z^{(a_i)}_r(\omega)= \z_{(r-V_{a_i}(\omega))+}(\omega^{[s_i]}).$$
The function $r\mapsto Z^{(a_i)}_r(\omega)$ then satisfies the properties stated in (i). Moreover it is immediate 
that the set of discontinuity times of this function is
$$\{V_{a_i}+r:r\in D(\omega^{[s_i]})\}= D^{(a_i)}(\omega)$$
where the last equality is again a consequence of our definitions. Futhermore, if $r\in D^{(a_i)}(\omega)$
and if $j$ is an index such that $a_j\in \check\cc^{(a_i)}_r$, the fact that
$\cc^{(a_i)}_{r-}$ is the closure of the union ${\cc}^{(a_i)}_r \cup \check\cc^{(a_i)}_r$ implies that
$$Z^{(a_i),\ve}_{r-}(\omega)= Z^{(a_i),\ve}_{r}(\omega)+ Z^{(a_j),\ve}_{r}(\omega)$$
and by passing to the limit $\ve\to 0$,
$$Z^{(a_i)}_{r-}(\omega)= Z^{(a_i)}_{r}(\omega)+ Z^{(a_j)}_{r}(\omega).$$ Finally property (d)
gives $Z^{(a_i)}_{r}(\omega)\not = Z^{(a_j)}_{r}(\omega)$.

The preceding discussion shows that properties (i) and (iii) of the proposition hold if
we restrict our attention to points in the dense sequence $a_1,a_2,\ldots$. However, it readily follows from our
definitions that we have $\cc^{(a)}_r=\cc^{(a_i)}_r$, and thus also $Z^{(a),\ve}_r=Z^{(a_i),\ve}_r$ as soon as $r<\min_{c\in\llbracket a,a_i\rrbracket}V_c$.
We infer that we can define $(Z^{(a)}_r)_{r\in [0,V_a)}$ in a unique way so that
$$Z^{(a)}_r=Z^{(a_i)}_r\;,\hbox{ for every } r<\min_{c\in\llbracket a,a_i\rrbracket}V_c\;,\hbox{ for every }i\geq 1.$$
It is then a simple matter to verify that assertions (i) and (iii) hold in the stated form, and assertion (ii) is also immediate.
\endproof

\section{The locally largest evolution}
\label{sec:local-larg}

We say that a point $a\in \t^\circ_\zeta$ is regular if
$$\bigcap_{t\in[0,V_a)} \cc^{(a)}_t = \{a\}.$$

\begin{proposition}
\label{loc-largest}
There exists $\N^{*,z}_0$ a.e. a unique point $a^\bullet$ of $\t^\circ_\zeta$ such that the
following two properties hold:
\begin{itemize}
\item[\rm(i)] We have
$Z^{(a^\bullet)}_t>|\Delta Z^{(a^\bullet)}_t|\;,\hbox{ for every }t\in D^{(a^\bullet)}$.
\item[\rm(ii)] The point $\ab$ is regular.
\end{itemize}
\end{proposition}

We will call $\ab$ the terminal point of the locally largest evolution.

\proof We first establish uniqueness. Suppose that $a_1$ and $a_2$ are two distinct points 
of $\t^\circ_\zeta$ that satisfy the properties stated in (i) and (ii). We notice that we must have
$$\min_{c\in \llbracket a_1,a_2\rrbracket} V_c < V_{a_1}\wedge V_{a_2}$$
because if the latter minimum is equal say  to $V_{a_1}$ the whole segment $\llbracket a_1,a_2\rrbracket$ is contained in 
$$\bigcap_{t\in[0,V_{a_1})} \cc^{(a_1)}_t $$
contradicting the regularity of $a_1$. Set $r=\min_{c\in \llbracket a_1,a_2\rrbracket} V_c $. By definition, we have then $r\in D^{(a_1)}\cap D^{(a_2)}$
and by property (iii) in Proposition \ref{pro-exit}, we get
$$Z^{(a_1)}_{r-}=Z^{(a_2)}_{r-}=Z^{(a_1)}_r + Z^{(a_2)}_r$$
and thus $|\Delta Z^{(a_1)}_r|=Z^{(a_2)}_r$, $|\Delta Z^{(a_2)}_r|=Z^{(a_1)}_r$. This shows that property (i) in Proposition \ref{loc-largest}
cannot hold simultaneously for $a_1$ and for $a_2$. 

Let us turn to existence. We let $t_\infty$ be the supremum of the set of all reals $t\geq 0$ such that there exists $a\in \t^\circ_\zeta$ with 
$V_a\geq t$ and $Z^{(a)}_s>|\Delta Z^{(a)}_s|$ for every $s<t$. By the definition of $t_\infty$, we can then find a 
nondecreasing sequence $(t_n)_{n\geq 1}$ in $\R_+$ and a corresponding sequence $(a_n)_{n\geq 1}$
in $\t^\circ_\zeta$ such that $V_{a_n}\geq t_n$, $Z^{(a_n)}_s>|\Delta Z^{(a_n)}_s|$ for every $s<t_n$, and $t_n\uparrow t_\infty$ as $n\uparrow \infty$. 
By compactness,we may assume that $a_n\la a_\infty\in\t_\zeta$ as $n\uparrow \infty$, and it is clear that $a_\infty\in\t^\circ_\zeta$
(the case $V_{a_\infty}=0$ is excluded since it would imply that $t_\infty=0$). 

Using property (ii) in Proposition \ref{pro-exit}, we have, for every $n$,
$$Z^{(a_\infty)}_s>|\Delta Z^{(a_\infty)}_s|\;,\quad\hbox{for every }s<t_n\wedge \min_{c\in \llbracket a_n,a_\infty\rrbracket} V_c.$$
Since
$$\min_{c\in \llbracket a_n,a_\infty\rrbracket} V_c \build{\la}_{n\to\infty}^{} V_{a_\infty}\geq t_\infty$$
(the last inequality because $V_{a_n}\geq t_n$ for every $n$), it follows that
$$Z^{(a_\infty)}_s>|\Delta Z^{(a_\infty)}_s|\;,\quad\hbox{for every }s<t_\infty.$$
We next claim that $t_\infty=V_{a_\infty}$. If not the case, we have $t_\infty<V_{a_\infty}$  and then either $t_\infty$
is a continuity time of $Z^{(a_\infty)}$, which immediately gives a contradiction with the
definition of $t_\infty$, or $t_\infty$
is a discontinuity time of $Z^{(a_\infty)}$, and by taking $a=a_\infty$ or $a\in \check\cc^{(a_\infty)}_{t_\infty}$ we
again get a contradiction with the definition of $t_\infty$. 

It remains to verify that $a_\infty$ is regular. 
If $a_\infty$ is not regular, then we can choose 
$$b\in \bigcap_{r\in[0,V_{a_\infty})} \cc^{(a_\infty)}_r$$
with $V_b>V_{a_\infty}=t_\infty$. If $t_\infty$ is a continuity point of $r\mapsto Z^{(b)}_r$,
we get a contradiction with the definition of $t_\infty$. If $t_\infty$ is a discontinuity point of $r\mapsto Z^{(b)}_r$,
then by taking $a=b$ or $a\in \check\cc^{(b)}_{t_\infty}$ we again get a contradiction. This completes the proof. \endproof

Let $u\in(0,V_{\ab})$. For future use, we note that, for every $a\in \t^\circ_\zeta$, we have
$a\in \cc^{(\ab)}_u$ if and only if $V_a>u$ and $Z^{(a)}_t >|\Delta Z^{(a)}_t|$ for every $t\leq u$. The ``only if'' part is trivial.
Conversely, assuming that $V_a>u$ and $Z^{(a)}_t >|\Delta Z^{(a)}_t|$ for every $t\leq u$, the property 
$a\notin \cc^{(\ab)}_u$ would lead to a contradiction by using Proposition \ref{pro-exit} (iii) with $r=\min_{c\in\llbracket \ab,a\rrbracket} V_c \leq u$.

\section{The law of the locally largest evolution}
\label{sec:law-loc}

Our next goal is to compute the distribution of $(Z^{(\ab)}_t)_{0\leq t < V_{\ab}}$ under $\N^{*,z}_0$. 

\begin{proposition}
\label{law-loc-largest}
The process $(Z^{(\ab)}_t)_{0\leq t < V_{\ab}}$ is distributed under $\N^{*,z}_0$ as $( X_t)_{0\leq t< H_0}$, where
$( X_t)_{t\geq 0}$ is the
self-similar Markov process  with index $\frac{1}{2}$ starting from $z$, which can be represented as
$$ X_t=z\,\exp(\xi(\chi(z^{-1/2}t))),$$
where $(\xi(s))_{s\geq 0}$ is the L\'evy process with only negative jumps whose Laplace exponent $\psi$ is given by formula 
\eqref{formula-psi} and
$(\chi(t))_{t\geq 0}$ is the time change defined in \eqref{formula-chi},
and $ H_0=\inf\{t\geq 0:  X_t=0\}$.
\end{proposition}

\proof
We fix $u>0$ and 
consider a bounded measurable function $F$ on the Skorokhod space
$\mathbb{D}([0,u],\R)$. We observe that
$$F\big((Z^{(\ab)}_t)_{0\leq t\leq u}\big)\,\mathbf{1}_{\{V_{\ab}>u\}}
= \int \vol(\dd a)\,F\big((Z^{(a)}_t)_{0\leq t\leq u}\big)\,\mathbf{1}_{\{a\in \cc^{(\ab)}_u\}}\,\frac{1}{\vol(\cc^{(a)}_u)},$$
simply because if $a\in \cc^{(\ab)}_u$ we have $Z^{(a)}_t=Z^{(\ab)}_t$ for $0\leq t\leq u$ and $\cc^{(a)}_u=\cc^{(\ab)}_u$.
From the definition of $\vol(\cdot)$ and a previous observation, the right-hand side can also be written as
$$\int_0^\sigma \dd s\,F\big((Z^{(s)}_t)_{0\leq t\leq u}\big)\,\mathbf{1}_{\{\wh W_s>u;\,Z^{(s)}_t>|\Delta Z^{(s)}_t|,\forall t\leq u\}}\,\frac{1}{\vol(\cc^{(p_\zeta(s))}_u)},$$
where we have written $Z^{(s)}_t=Z^{(a)}_t$ if $a=p_\zeta(s)$ to simplify notation. 

The preceding considerations show that
\begin{align*}
&\N^{*,z}_0\Big(F\big((Z^{(\ab)}_t)_{0\leq t\leq u}\big)\,\mathbf{1}_{\{V_{\ab}>u\}}\Big)\\
&=\N^{*,z}_0\Big(\int_0^\sigma \dd s\,F\big((Z^{(s)}_t)_{0\leq t\leq u}\big)\,
\mathbf{1}_{\{\wh W_s>u;\,Z^{(s)}_t>|\Delta Z^{(s)}_t|,\forall t\leq u\}}\,\frac{1}{\vol(\cc^{(p_\zeta(s))}_u)}\Big)\\
&= z^2\,\N_0\Big(F\big(\z_{(L_z+t)+})_{0\leq t\leq u}\big)\,
\mathbf{1}_{\{L_z<-u;\,\z_{(L_z+t)+}>|\Delta\z_{L_z+t}|,\forall t\leq u\}}\, \frac{1}{\vol(\cc_{L_z+u})}
\,\Big|\,L_z<0\Big),
\end{align*}
where for $r<0$, we use (under $\N_0$) the notation $\cc_r$ for the connected component of $\{a\in\t_\zeta: V_a > r\}$
containing the ``root'' $p_\zeta(0)$. The second equality of the last display is a consequence of Proposition \ref{repre}
and the way the functions $Z^{(a)}_t$ were constructed in Section \ref{sec:bdry-size}.

In the terminology of \cite{ALG}, $\cc_{L_z+u}$ is (up to a set of zero volume) the union of the subsets of $\t_\zeta$
corresponding to the excursions above the minimum that start at a level greater than $L_z+u$. Using Theorem \ref{mainALG} and 
\cite[Proposition 31]{ALG}, we obtain that the conditional distribution of $\vol(\cc_{L_z+u})$ under $\N_0(\cdot\mid L_z<0)$ and
knowing $(\z_r)_{r<0}$  is the law of
$$\sum_{i=1}^\infty |\Delta\z_{r_i}|^2 \,\nu_i$$
where $r_1,r_2,\ldots$ is an enumeration of the jumps of $\z$ on $(L_z+u,0)$, and the random variables $\nu_1,\nu_2,\ldots$
are independent and distributed according to the density
$$\frac{1}{\sqrt{2\pi}}\,x^{-5/2}\,\exp(-\frac{1}{2x})\,\mathbf{1}_{\{x>0\}}.$$
Writing $\E^{(\nu)}[\cdot]$ for the expectation with respect to the variables $\nu_1,\nu_2,\ldots$, we can thus also write
\begin{align}
\label{key1}
&\N^{*,z}_0\Big(F\big((Z^{(\ab)}_t)_{0\leq t\leq u}\big)\,\mathbf{1}_{\{V_{\ab}>u\}}\Big)\nonumber\\
&\quad= z^2\,\N_0\Big(F\big((\z_{L_z+t})_{0\leq t\leq u}\big)\,
\mathbf{1}_{\{L_z<-u;\,\z_{(L_z+t)+}>|\Delta\z_{L_z+t}|,\forall t\leq u\}}\, 
\E^{(\nu)}\Big[ \frac{1}{\sum  |\Delta\z_{r_i}|^2 \,\nu_i}\Big]\,\Big|\,L_z<0\Big),
\end{align}
and the right-hand side is an integral under $\N_0$ of a 
quantity depending only on the exit measure process $(\z_r)_{r<0}$.

\medskip
Thanks to Proposition \ref{law-exit-p}, we can replace the right-hand side of \eqref{key1} by
\begin{equation}
\label{key2}
z^2\,\E\Big[F\big((X^\circ_t)_{0\leq t\leq u})\,\mathbf{1}_{\{H^\circ_0>u;\,X^\circ_t>|\Delta X^\circ_t|,\forall t\leq u\}}\, 
\E^{(\nu)}\Big[ \frac{1}{\sum  |\Delta X^\circ_{s_i}|^2 \,\nu_i}\Big]\Big],
\end{equation}
where $s_1,s_2,\ldots$ is an enumeration of the jump times of $X^\circ$ over $[u,H^\circ_0)$.
Then, by the Markov property and the self-similarity of $X^\circ$, the conditional expectation of
the quantity
$$\E^{(\nu)}\Big[ \frac{1}{\sum  |\Delta X^\circ_{s_i}|^2 \,\nu_i}\Big]$$
given $(X^\circ_t)_{0\leq t\leq u}$ is $C/(X^\circ_u)^2$, for some constant $C>0$ (the cases $C=0$ and $C=\infty$
are excluded since the preceding equalities would give an absurd statement).
So the quantity \eqref{key2} is also equal to $C$ times
\begin{equation}
\label{key3}
z^2\,\E\Big[F\big((X^\circ_t)_{0\leq t\leq u}\big)\,\mathbf{1}_{\{H^\circ_0>u;\,X^\circ_t>|\Delta X^\circ_t|,\forall t\leq u\}}\, 
(X^\circ_u)^{-2}\Big],
\end{equation}

We will rewrite this quantity in a different form. In the remaining part of the proof, we take $z=1$ for the sake of simplicity (of course the
self-similarity of $X^\circ$ will then allow us to get a similar result for an arbitrary value of $z$).
Using the representation in Proposition \ref{law-exit-p}, we obtain that the quantity \eqref{key3}
is equal for $z=1$ to
\begin{equation}
\label{key5}
\E\Big[F\Big(\big(\exp(\xi^\circ(\chi^\circ(t)))\big)_{0\leq t\leq u}\Big)\,\mathbf{1}_{\{\chi^\circ(u)<\infty\}}\,
\mathbf{1}_{\{\Delta \xi^\circ(s)>-\log 2,\,\forall s\in [0,\chi^\circ(u)]\}}
\,\exp(-2\xi^\circ(\chi^\circ(u)))\Big].
\end{equation}

\begin{lemma}
\label{chang-pro}
For every $v\geq 0$, set 
$$M_v= \mathbf{1}_{\{\Delta \xi^\circ(s)>-\log 2,\forall s\in [0,v]\}}
\,\exp(-2\xi^\circ(v)).$$
Then $(M_v)_{v\geq 0}$ is a martingale with respect to the canonical filtration of the process $\xi^\circ$. Let 
$ \xi$ be as in Proposition \ref{law-loc-largest}.
Then, for every fixed $v>0$ the process 
$(\xi^\circ(t))_{0\leq t\leq v}$ is distributed under the probability measure $M_v\cdot \P$ as $(\xi(t))_{0\leq t\leq v}$ under $\P$. 
\end{lemma}

\proof To simplify notation, we write $\alpha=\sqrt{3/2\pi}$. Thanks to the properties of L\'evy processes, 
in order to verify that  $(M_v)_{v\geq 0}$ is a martingale, it suffices to prove that $\E[M_v]=1$ for every $v>0$. 
It is convenient to set
$$\xi''(t)=\sum_{0\leq s \leq t} \Delta \xi^\circ(s)\, \mathbf{1}_{\{\Delta\xi^\circ(s)\leq -\log 2\}},$$
so that we can write $\xi^\circ(t)=\xi'(t) + \xi''(t)$, where $\xi'$ and $\xi''$ are two independent L\'evy processes.
The Laplace exponent of $\xi''$ is
$$\psi''(\lambda)=\alpha \int_{-\infty}^{-\log 2} (e^{\lambda y}-1)\,e^{y/2}\,(1-e^y)^{-5/2}\,\dd y,$$
and the Laplace exponent of $\xi'$ is
\begin{align}
\label{key4}
\psi'(\lambda)&=\psi^\circ(\lambda)-\psi''(\lambda)\nonumber\\
&=\alpha\Bigg( \int_{-\log 2}^0 (e^{\lambda y}-1-\lambda(e^y-1))\,e^{y/2}\,(1-e^y)^{-5/2}\,\dd y
-\lambda \int_{-\infty}^{-\log2} (e^y-1)\,e^{y/2}\,(1-e^y)^{-5/2}\,\dd y\Bigg)\nonumber\\
&=\alpha\Bigg( \int_{-\log 2}^0 (e^{\lambda y}-1-\lambda(e^y-1))\,e^{y/2}\,(1-e^y)^{-5/2}\,\dd y + 2\lambda\Bigg)
\end{align}
using the simple calculation
$$\int_{-\infty}^{-\log 2} e^{y/2}\,(1-e^y)^{-3/2}\,\dd y =\int_0^{1/2} (1-x)^{-3/2}\,\frac{\dd x}{\sqrt{x}}= 2.$$
Note that $\xi'$ has bounded jumps and therefore exponential moments of any order, so that 
$\psi'(\lambda)$ makes sense for every $\lambda \in\R$ and not only $\lambda \geq 0$.

We have then 
$$\E[M_v]=\P(\xi''(v)=0)\,\E[e^{-2\xi'(v)}].$$
On one hand, $\E[e^{-2\xi'(v)}]=\exp(\alpha K v)$, where
\begin{align*}
K&=-4 + \int_{-\log 2}^0 (e^{-2 y}-1+2(e^y-1))\,e^{y/2}\,(1-e^y)^{-5/2}\,\dd y\\
&=-4+ \int_{1/2}^1 (x^{-2}-1+2(x-1))\,(1-x)^{-5/2}\,\frac{\dd x}{\sqrt{x}}\\
&=-4 + 2\int_0^{1/2} x^{-3/2}\,(1-x)^{-1/2}\,\dd x + \int_{1/2}^1 x^{-5/2}\,(1-x)^{-1/2}\,\dd x \\
&=\frac{8}{3}
\end{align*}
and on the other hand
$$\P(\xi''(v)=0)=\exp\Big(-\alpha v\int_{-\infty}^{-\log2} e^{y/2}\,(1-e^y)^{-5/2} \dd y\Big)=\exp(-\frac{8}{3}\alpha v).$$
By combining the last two displays, we get the desired result $\E[M_v]=1$. 

Then, let us fix $v>0$. It is straightforward to verify that the properties of stationarity and independence of the 
increments of $\xi^\circ$ are preserved under the probability measure $M_v\cdot \P$, so that $(\xi^\circ(t))_{0\leq t\leq v}$
remains a L\'evy process under this probability measure. To evaluate the Laplace
exponent of this L\'evy process, we write
$$
\E[M_v\,e^{\lambda \xi^\circ(v)}]=\exp(-\frac{8}{3}\alpha v)\,\E[e^{(\lambda-2)\xi'(v))}]=\exp((\psi'(\lambda-2)-\frac{8}{3}\alpha)v)
=\exp(v\psi(\lambda)),
$$
where $\psi$ is as in \eqref{formula-psi}. The last equality follows from formula \eqref{key4} for $\psi'$ and simple calculations left to the reader. This completes the proof
of the lemma. \endproof

Let us come back to \eqref{key5}. We first observe that, for every $v>0$,
\begin{align}
\label{key6}
\E\Big[F\Big(\big(\exp(\xi^\circ(\chi^\circ (t)))\big)_{0\leq t\leq u}\Big)\,\mathbf{1}_{\{\chi^\circ(u)\leq v\}}\,M_{\chi^\circ(u)}\Big]
&=\E\Big[F\Big(\big(\exp(\xi^\circ(\chi^\circ (t)))\big)_{0\leq t\leq u}\Big)\,\mathbf{1}_{\{\chi^\circ(u)\leq v\}}\,M_{v}\Big]\nonumber\\
&=\E\Big[F\Big(\big(\exp(\xi(\chi(t)))\big)_{0\leq t\leq u}\Big)\,\mathbf{1}_{\{\chi(u)\leq v\}}\Big],
\end{align}
where $\chi(t)$ is as in formula \eqref{formula-chi}.
As $v\to\infty$, the right-hand side of \eqref{key6} converges to the quantity $\E[F(\exp(\xi(\chi(t))))_{0\leq t\leq u})\,\mathbf{1}_{\{\chi(u)<\infty\}}]$.
Similarly the left-hand side of \eqref{key6} converges to the quantity \eqref{key5}: Dominated convergence is easily justified by noting that
$$\E[M_{\chi^\circ(u)}\,\mathbf{1}_{\{\chi^\circ(u)<\infty\}}]\leq \liminf_{t\to\infty} \E[M_{\chi^\circ(u)\wedge t}\,\mathbf{1}_{\{\chi^\circ(u)\leq t\}}]$$
and $\E[M_{\chi^\circ(u)\wedge t}]=\E[M_t]=1$ by the optional stopping theorem.

Hence the quantity \eqref{key5} is equal to $\E[F\big((\exp(\xi(\chi(t))))_{0\leq t\leq u}\big)\,\mathbf{1}_{\{\chi(u)<\infty\}}]$.
We conclude that
$$\N^{*,1}_0\Big(F\big((Z^{(\ab)}_t)_{0\leq t\leq u}\big)\,\mathbf{1}_{\{V_{\ab}>u\}}\Big)= C\,\E\Big[F\big((\exp(\xi(\chi(t))))_{0\leq t\leq u}\big)\,\mathbf{1}_{\{\chi(u)<\infty\}}\Big].$$
At this stage, we can take $F=1$ and let $u$ tend to $0$, and we find that $C=1$. This gives the statement of Proposition \ref{law-loc-largest} for $z=1$, and
it is easily extended by self-similarity. \hfill$\square$

\section{Excursions from the locally largest evolution}
\label{sec:excu-lle}

If $\omega\in\S_0$ satisfies the regularity properties stated in Section \ref{comp-level}, we can define the excursions above the minimum 
$\omega^r$, $r\in D$ in the way
described in Section \ref{excu-mini}. 

Let us now argue under $\N^{*,z}_0(\dd \omega)$ and for $u\in[0,\sigma]$, consider the re-rooted snake trajectory $\omega^{[u]}$. Let $a=p_\zeta(u)$
and recall the set $D^{(a)}$ corresponding to the discontinuity times of 
$(Z^{(a)}_r)_{r\in [0,V_a)}$. As we already noticed in the proof of Proposition \ref{pro-exit}, we have
$D^{(a)}(\omega)=\{V_a+r: r\in D(\omega^{[u]})\}$. If $r\in D^{(a)}(\omega)$ we can thus associate 
with $r-V_a\in D(\omega^{[u]})$ an excursion of $\omega^{[u]}$ above the minimum, which we
denote by ${\omega}^{a,r}$ (it is easy to see that this excursion only depends on $a$ and not on
$u$ such that $a=p_\zeta(u)$). We already noticed that, if $a,a'\in \t^\circ_\zeta$ are such that $V_a\wedge V_{a'}>u$ and
$\cc^{(a)}_u=\cc^{(a')}_u$, we have
$D^{(a)}\cap [0,u]=D^{(a')}\cap[0,u]$, and it is also true that $\omega^{a,r}=\omega^{a',r}$
for $r\in D^{(a)}\cap [0,u]$. 

We will now apply the preceding considerations to $a=\ab$. We write $D^{(\ab)}=\{r_1,r_2,\ldots\}$, where 
$|\Delta Z^{(\ab)}_{r_1}|> |\Delta Z^{(\ab)}_{r_2}|>\cdots$.

\begin{proposition}
\label{excu-loc-largest}
Under $\N^{*,z}_0$, conditionally on $(Z^{(\ab)}_r)_{0\leq r< V_\ab}$, the excursions $\omega^{\ab,r_i}$, $i=1,2,\ldots$
are independent, and for every fixed $i\geq 1$ the conditional distribution of $\omega^{\ab,r_i}$ is $\N^{*,|\Delta Z^{(\ab)}_{r_i}|}_0$.
\end{proposition}

\proof We proceed in a way similar to the one used above to determine the law of $(Z^{(\ab)}_r,0\leq r\leq V_\ab)$.
We fix $u>0$ and consider a bounded measurable function $F$ on the Skorokhod space
$\D([0,u],\R_+)$, and a bounded measurable function $H$ on $\R_+\times\S_0$ such that $H=1$ on $\R_+\times (\S_0\backslash \S^{(\delta)}_0)$
for some $\delta>0$. The latter condition ensures that $H(r,\omega^{\ab,r})=1$ except for finitely many values 
of $r\in D$. Then,
\begin{align}
\label{loc-lar-excu}
&\N^{*,z}_0\Big(F\big((Z^{(\ab)}_t)_{0\leq t\leq u}\big)\,\mathbf{1}_{\{V_{\ab}>u\}}\,\prod_{r\in D^{(\ab)}\cap[0,u]} H(r,\omega^{\ab,r})\Big)\nonumber\\
&\quad=\N^{*,z}_0\Big(\int_0^\sigma \dd s\,F\big((Z^{(s)}_t)_{0\leq t\leq u}\big)\,
\mathbf{1}_{\{\wh W_s>u;\,Z^{(s)}_t>|\Delta Z^{(s)}_t|,\forall t\leq u\}}\,\nonumber\\
&\qquad\qquad\quad\times\frac{1}{\vol(\cc^{(p_\zeta(s))}_u)}
\,\prod_{r\in D^{(p_\zeta(s))}\cap[0,u]} H(r,\omega^{p_\zeta(s),r})\Big)\nonumber\\
&\quad= z^2\,\N_0\Big(F\big((\z_{(L_z+t)+})_{0\leq t\leq u}\big)\,\mathbf{1}_{\{L_z<-u;\,\z_{(L_z+t)+}>|\Delta\z_{L_z+t}|,\forall t\leq u\}}\,\,\nonumber\\
&\qquad\qquad\quad\times \frac{1}{\vol(\cc_{L_z+u})}
\prod_{r\in D\cap[L_z,L_z+u]} H(r-L_z,\omega^{r})\Big).
\end{align}
We used the remarks preceding the statement of the proposition in the first equality, and Proposition \ref{repre} in the second one. At this stage, we use 
Theorem \ref{mainALG}, which shows that conditionally 
on the exit measure process $(\z_r)_{r<0}$ (whose set of discontinuities is $D$) the excursions $\omega_r$, $r\in D$, are
independent and the conditional distribution of $\omega^r$ is $\N^{*,|\Delta \z_r|}_0$. Since the quantity $\vol(\cc_{L_z+u})$
only depends on the excursions $\omega^r$ with $r>L_z+u$, we can rewrite the last line of the preceding display as
\begin{align*}
&z^2\,\N_0\Big(F\big((\z_{(L_z+t)+})_{0\leq t\leq u}\big)\,\mathbf{1}_{\{L_z<-u;\,\z_{(L_z+t)+}>|\Delta\z_{L_z+t}|,\forall t\leq u\}}\,\,\\
&\qquad\quad\times\frac{1}{\vol(\cc_{L_z+u})}
\prod_{r\in D\cap[L_z,L_z+u]} \N^{*,|\Delta \z_r|}_0(H(r-L_z,\cdot))\Big)
\end{align*}
but then, we can re-use the same arguments ``backwards'' to see that the latter quantity is also equal to
\begin{equation}
\label{loc-lar-excu2}
\N^{*,z}_0\Big(F\big((Z^{(\ab)}_t)_{0\leq t\leq u}\big)\,\mathbf{1}_{\{V_{\ab}>u\}}\,\prod_{r\in D^{(\ab)}\cap[0,u]}  \N^{*,|\Delta Z^{(\ab)}_r|}_0(H(r,\cdot))\Big).
\end{equation}
The quantity \eqref{loc-lar-excu2} is thus equal to the left-hand-side of \eqref{loc-lar-excu}. Simple arguments show that this also
implies that, for any bounded measurable function $G$ on the appropriate space of c\`adl\`ag paths,
$$
\N^{*,z}_0\Big(G\big((Z^{(\ab)}_t)_{0\leq t<V_{\ab}}\big)\!\prod_{r\in D^{(\ab)}} H(r,\omega^{\ab,r})\Big)
=\N^{*,z}_0\Big(G\big((Z^{(\ab)}_t)_{0\leq t<V_{\ab}}\big)\!\prod_{r\in D^{(\ab)}} \N^{*,|\Delta Z^{(\ab)}_r|}_0(H(r,\cdot))\Big).
$$
This gives the statement of the proposition.
 \endproof
 
 We will call $\omega^{\ab,r_i}$, $i=1,2,\ldots$ the excursions of $\omega$ from the locally largest evolution. The number 
 $r_i$ is called the starting level of $\omega^{\ab,r_i}$.
As previously, these excursions will always be listed in decreasing order of their ``boundary sizes'' $|\Delta Z^{(\ab)}_{r_i}|$.

\section{The growth-fragmentation process}
\label{gro-frag}

In this section, we argue again under $\N^{*,z}_0(\dd \omega)$. To simplify notation, we will write $\omega^{(1)},\omega^{(2)},\ldots$
for the excursions of $\omega$ from the locally largest evolution ($\omega^{(i)}=\omega^{\ab,r_i}$ in the notation 
of the previous section), but we keep the notation $r_1,r_2,\ldots$ for the respective starting levels of these excursions.

Next, for every $i\geq 1$, since the conditional distribution of $\omega^{(i)}$ knowing $(Z^{(\ab)}_r)_{0\leq r< V_\ab}$ 
is $\N^{*,|\Delta Z^{(\ab)}_{r_i}|}_0$, we can also define a point $a^\bullet_{(i)}$ as the terminal point 
of the locally largest evolution in $\omega^{(i)}$, and the excursions $\omega^{(i,1)},\omega^{(i,2)},\ldots$ from the locally largest
evolution in $\omega^{(i)}$ (ranked as explained at the end of the previous section). We write $r_{i,j}$ for the starting level of $\omega^{(i,j)}$.

Obviously we can continue the construction by induction. Assuming that we have defined $\omega^{(i_1,\ldots,i_k)}$,
we let $a^\bullet_{(i_1,\ldots,i_k)}$ be the terminal point 
of the locally largest evolution in $\omega^{(i_1,\ldots,i_k)}$, and we denote the excursions from 
the locally largest
evolution in $\omega^{(i_1,\ldots,i_k)}$ by $\omega^{(i_1,\ldots,i_k,1)},\omega^{(i_1,\ldots,i_k,2)},\ldots$. For every $j\geq 1$,
we let $r_{i_1,\ldots,i_k,j}$ be the starting level of $\omega^{(i_1,\ldots,i_k,j)}$.

We also set, for every $(i_1,\ldots,i_k)$,
$$h_{i_1,\ldots,i_k}= r_{i_1}+r_{i_1,i_2}+ \cdots + r_{i_1,\ldots,i_k},$$
and we let $\beta_{i_1,\ldots,i_k}$ be the label of $a^\bullet_{(i_1,\ldots,i_k)}$ (in $\omega^{(i_1,\ldots,i_k)}$).

Let $r\in[ h_{i_1,\ldots,i_k}, h_{i_1,\ldots,i_k}+\beta_{i_1,\ldots,i_k})$, and consider the connected component $\cc^{(a^\bullet_{(i_1,\ldots,i_k)})}_{r-h_{i_1,\ldots,i_k}}(\omega^{(i_1,\ldots,i_k)})$
(as in Section \ref{comp-level}, this is the connected component of $\{a\in \t_\zeta(\omega^{(i_1,\ldots,i_k)}):V_a(\omega^{(i_1,\ldots,i_k)})>r-h_{i_1,\ldots,i_k}\}$ 
that contains $a^\bullet_{(i_1,\ldots,i_k)}$). As explained in Section \ref{excu-mini}, this connected component corresponds (via a volume-preserving
isometry) to a connected component of $\{a\in \t_\zeta(\omega^{(i_1,\ldots,i_{k-1})}):V_a(\omega^{(i_1,\ldots,i_{k-1})})>r-h_{i_1,\ldots,i_{k-1}}\}$
and inductively to a connected component of $\{a\in\t_\zeta(\omega): V_a(\omega)>r\}$. The latter component is denoted by 
$\mathcal{D}^{(i_1,\ldots,i_k)}_r$. Recall that this definition makes sense only if $r\in[ h_{i_1,\ldots,i_k}, h_{i_1,\ldots,i_k}+\beta_{i_1,\ldots,i_k})$
(otherwise we may take $\mathcal{D}^{(i_1,\ldots,i_k)}_r=\varnothing$).

We set
$$\mathcal{U}=\bigcup_{k=0}^\infty \N^k$$
with the convention $\N^0=\{\varnothing\}$.
We define $h_\varnothing=0$, $\beta_\varnothing=V_\ab(\omega)$ and $\mathcal{D}^\varnothing_r=\cc^{(\ab)}_r$ if $0\leq r<V_\ab(\omega)$.

\begin{lemma}
\label{all-cc}
Let $r\geq 0$.
The sets $\mathcal{D}^{(i_1,\ldots,i_k)}_r$, for all $(i_1,\ldots,i_k)\in \mathcal{U}$ such that $h_{i_1,\ldots,i_k}\leq r< h_{i_1,\ldots,i_k}+\beta_{i_1,\ldots,i_k}$,
are exactly the connected components of $\{a\in \t_\zeta: V_a>r\}$.
\end{lemma}

\proof
We already know that any of the sets $\mathcal{D}^{(i_1,\ldots,i_k)}_r$, $(i_1,\ldots,i_k)\in \mathcal{U}$, is a connected component of 
$\{a\in \t_\zeta: V_a>r\}$, and we need to show that any connected component is of this type.
Let $\cc$ be a connected component of 
$\{a\in \t_\zeta: V_a>r\}$, and choose any $a\in \cc$. The process $(Z^{(a)}_t)_{0\leq t\leq r}$ has only finitely many jump times $s\in[0,r]$
such that $|\Delta Z^{(a)}_s|>Z^{(a)}_s$, and we denote these jump times by $0<t_1<t_2<\cdots<t_k\leq r$, where $k\geq 0$.

\smallskip
$\bullet$ If $k=0$, this means that $|\Delta Z^{(a)}_s|<Z^{(a)}_s$ for every $s\in[0,r]$, and we have already seen that
this implies that $a\in \cc^{(a^\bullet)}_r$. We have thus $\cc=\cc^{(a^\bullet)}_r=\mathcal{D}^\varnothing_r$ in that case.

\smallskip
$\bullet$ Suppose that $k\geq 1$. We have $Z^{(a)}_s=Z^{(\ab)}_s$ if and only if $0\leq s<t_1$. In particular $a\in\cc^{(a^\bullet)}_s$
if and only if $0\leq s<t_1$, so that $a$ belongs to $a\in\cc^{(a^\bullet)}_{t_1-}$, which is the closure of $\cc^{(a^\bullet)}_{t_1}\cup \check\cc^{(a^\bullet)}_{t_1}$.
Since points of the boundary of $\cc^{(a^\bullet)}_{t_1}\cup \check\cc^{(a^\bullet)}_{t_1}$ have label $t_1$ whereas $V_a>t\geq t_1$, it follows that
$a\in \check\cc^{(a^\bullet)}_{t_1}$, and $\cc\subset\check\cc^{(a^\bullet)}_{t_1}$. Furthermore $t_1$ is a jump time of $Z^{(\ab)}$, so that we have
$t_1=r_{i_1}$ for some $i_1\geq1$. As explained in Section \ref{excu-mini}, $\check\cc^{(a^\bullet)}_{t_1}$ is identified to $\t^\circ_\zeta(\omega^{(i_1)})$,
and through this identification $\cc$ is identified to a connected component $\cc'$ of $\{b\in \t^\circ_\zeta(\omega^{(i_1)}):V_b(\omega^{(i_1)})>r-r_{i_1}\}$
and $a$ is identified to a point $a'$ of $\cc'$. We have then
$$(Z^{(a')}_t(\omega^{(i_1)}),0\leq t\leq r-t_1)= (Z^{(a)}_{r_1+t},0\leq t\leq r-t_1).$$
In particular, if $k=1$, there are no jump times $s\in[0,r-t_1]$ such that $|\Delta Z^{(a')}_s(\omega^{(i_1)})|>Z^{(a')}_s(\omega^{(i_1)})$, and
we conclude that $\cc'=\cc^{(a^\bullet_{(i_1)})}_{r-r_{i_1}}(\omega^{(i_1)})$, which means that $\cc=\mathcal{D}^{(i_1)}_r$.

\smallskip
$\bullet$ Suppose $k\geq 2$. Then $t_2-t_1$ is the first jump time 
of $(Z^{(a')}_t(\omega^{(i_1)}))_{0\leq t\leq r-t_1}$ such that $|\Delta Z^{(a')}_s(\omega^{(i_1)})|>Z^{(a')}_s(\omega^{(i_1)})$, 
we have $Z^{(a^\bullet_{(i_1)})}_s(\omega^{(i_1)})=Z^{(a')}_s(\omega^{(i_1)})$ for $0\leq s<t_2-t_1$, and there exists
$i_2\geq 1$ such that $t_2-t_1=r_{i_1,i_2}$. We have then $\cc'\subset \check\cc^{(a^\bullet_{(i_1)})}_{t_2-t_1}(\omega^{(i_1)})$.
It follows that $\cc'$ is identified to a connected component $\cc''$ of $\{b\in \t^\circ_\zeta(\omega^{(i_1,i_2)}):V_b(\omega^{(i_1,i_2)})>r-r_{i_1}-r_{i_1,i_2}\}$.
If $k=2$, we conclude as in the preceding step that $\cc''=\cc^{(a^\bullet_{(i_1,i_2)})}_{r-h_{i_1,i_2}}(\omega^{(i_1,i_2)})$, which means
that $\cc=\mathcal{D}^{(i_1,i_2)}_r$.

\smallskip
\noindent The proof is easily completed by induction, and we omit the details. \endproof

\noindent{\it Proof of Theorem \ref{main2}.} The first part of Theorem \ref{main2} (concerning the definition and approximation of
boundary sizes) is a consequence of Proposition \ref{pro-exit}, which in fact gives a stronger 
result. So we only need to prove the second part of the statement. If $\cc$ is a connected component of $\{a\in\t_\zeta: V_a>r\}$, we write
$Z_{(\cc)}=Z^{(a)}_r$ where $a$ is any point of $\cc$ (this does not depend on the choice of $a$). To simplify notation, 
for every $(i_1,\ldots,i_k)\in\mathcal{U}$ and every $j\geq 1$, we write $\Delta_{(i_1,\ldots,i_k,j)}$ for the jump at
time $h_{i_1,\ldots,i_k,j}$ of the function $r\mapsto Z_{(\mathcal{D}^{(i_1,\ldots,i_k)}_r)}$ --- from our construction this is
also the jump at time $r_{i_1,\ldots,i_k,j}$ of the function $r\mapsto Z^{(a^\bullet_{(i_1,\ldots,i_k)})}_{r}(\omega^{(i_1,\ldots,i_k)})$.

From 
the preceding lemma, we get that $\mathbf{Y}(r)$ is obtained as the (reordered) collection of the 
quantities $Z_{(\mathcal{D}^{(i_1,\ldots,i_k)}_r)}$ for all $(i_1,\ldots,i_k)\in\mathcal{U}$ such that 
$h_{i_1,\ldots,i_k}\leq r< h_{i_1,\ldots,i_k}+\beta_{i_1,\ldots,i_k}$. 

We know that the process $(Z_{(\mathcal{D}^{\varnothing}_r)})_{0\leq r<\beta_\varnothing}=(Z^{(\ab)}_r)_{0\leq r< V_\ab}$ is distributed as
the self-similar Markov process $ X$ of Proposition \ref{law-loc-largest} started from $z$ and killed when it hits $0$. Thanks to
Proposition \ref{excu-loc-largest}, we then get that,
conditionally on $(Z_{(\mathcal{D}^{\varnothing}_r)})_{0\leq r<\beta_\varnothing}$, the excursions $\omega^{(i)}$, $i=1,2,\ldots$, are
independent and for every fixed $j$ the conditional distribution of $\omega^{(j)}$ is $\N^{*,|\Delta_{(j)}|}_0$. Consequently, under the
same conditioning, the processes $(Z_{(\mathcal{D}^{(i)}_{h_i+r})})_{0\leq r\leq \beta_i}$, $i=1,2,\ldots$, are independent
copies of $ X$ started respectively at $|\Delta_{(i)}|$, $i=1,2,\ldots$

We can continue by induction, using Proposition \ref{excu-loc-largest} at every step. We obtain that,
conditionally on the processes
$$\Big(Z_{(\mathcal{D}^{(i_1,\ldots,i_\ell)}_r)}\Big)_{ h_{(i_1,\ldots,i_\ell)}\leq r<h_{(i_1,\ldots,i_\ell)} + \beta_{(i_1,\ldots,i_\ell)}},
\;0\leq \ell\leq k,(i_1,\ldots,i_\ell)\in\N^\ell,$$
the excursions $\omega^{(i_1,\ldots,i_k,j)}$, $(i_1,\ldots,i_k,j)\in\N^{k+1}$, are
independent and for every fixed $(i_1,\ldots,i_k,j)$ the conditional distribution of $\omega^{(i_1,\ldots,i_k,j)}$ is $\N^{*,|\Delta_{(i_1,\ldots,i_k,j)}|}_0$.
Hence, under the same conditioning, the processes 
$$\Big(Z_{(\mathcal{D}^{(i_1,\ldots,i_k,j)}_{h_{(i_1,\ldots,i_k,j)}+r})}\Big)_{0\leq r\leq \beta_{(i_1,\ldots,i_k,j)}},\ (i_1,\ldots,i_k,j)\in\N^{k+1}$$
are independent
copies of $ X$ started respectively at $|\Delta_{(i_1,\ldots,i_k,j)}|$.

From these observations,
we conclude that $(\mathbf{Y}(r))_{r\geq 0}$ is a growth fragmentation process whose 
Eve particle process is the self-similar Markov process $ X$ and with initial value $\mathbf{Y}(0)=(z,0,0,\ldots)$. \hfill$\square$

\section{Slicing the Brownian disk at heights}
\label{sec:slicing}

In this section, we prove Theorem \ref{main3}. We rely on the construction of 
the free Brownian disk $\D_z$ from a random snake trajectory distributed according
to $\N^{*,z}_0$. This construction is given in \cite{Disks}, to which 
we refer for additional details. Throughout this section, we argue under
$\N^{*,z}_0$, and the following statements hold $\N^{*,z}_0$ a.s.

The free Brownian disk $\D_z$ is a random geodesic compact metric space, which is constructed (under
$\N^{*,z}_0$) as a quotient space of $\t_\zeta$. The 
canonical projection, which is a continuous mapping from $\t_\zeta$ onto $\D_z$,
is denoted by $\Pi$. We note that the push forward of the volume measure $\mathrm{vol}$
on $\t_\zeta$ is the volume measure $\mathbf{V}$ on $\D_z$.

Recall 
the notation $H(x)$, for the ``height'' of $x\in\D_z$ (the distance from
$x$ to the boundary $\partial\D_z$). We will not need the details 
of the construction of $\D_z$, but we record the following two facts:
\begin{enumerate}
\item[(a)] If $a\in\t_\zeta$ and $x=\Pi(a)$, we have $H(x)=V_a$.
\item[(b)] For every $a,b\in \t_\zeta$ such that $\Pi(a)=\Pi(b)$, we have
$V_a=V_b=\displaystyle{\min_{c\in\llbracket a,b\rrbracket}} V_c$.
\end{enumerate}

The following lemma is an analog for the Brownian disk of Proposition
3.1 of \cite{Acta} for the Brownian map. The proof is similar, but we provide 
details because this result is the key to the derivation of Theorem \ref{main3}.

\begin{lemma}
\label{cactus}
Let $a,b\in \t_\zeta$ and let $(\gamma(t))_{0\leq t\leq T}$ be a continuous path in $\D_z$ such that
$\gamma(0)=\Pi(a)$ and $\gamma(T)=\Pi(b)$. Then 
$$\min_{0\leq t\leq T} H(\gamma(t))\leq \min_{c\in\llbracket a,b \rrbracket} V_c.$$
\end{lemma}

\proof We may assume that 
$$V_a\wedge V_b > \min_{c\in\llbracket a,b \rrbracket} V_c$$
since the result is trivial otherwise. Then we can find $c_0\in\rrbracket a,b\llbracket$ such that
$$V_{c_0}= \min_{c\in\llbracket a,b \rrbracket} V_c.$$
The points $a$ and $b$ are in different connected components of $\t_\zeta\backslash\{c_0\}$.
Let $\cc_1$ be the connected component of $\t_\zeta\backslash\{c_0\}$ that contains $a$,
and let $\cc_2=\t_\zeta\backslash \ov{\cc}_1$, so that $b\in \cc_2$. 
Set
$$t_0:=\inf\{t\in[0,T]: \gamma(t)\in \Pi(\ov{\cc}_2)\}.$$
Since $\Pi(\ov{\cc}_2)$ is closed, we have $\gamma(t_0)\in \Pi(\ov{\cc}_2)$. Furthermore, $t_0>0$
because otherwise this would mean that $\Pi(a)\in \Pi(\ov{\cc}_2)$, and thus there exists $a'\in \ov{\cc}_2$
such that $\Pi(a)=\Pi(a')$: Noting that $c_0\in\llbracket a,a'\rrbracket$, property (b) above would 
imply that $V_a\leq V_{c_0}$, which is a contradiction.

We can then choose a sequence $(s_n)_{n\geq 1}$ in $[0,t_0)$ such that $s_n\uparrow t_0$ as 
$n\uparrow\infty$. Since $\gamma({s_n})\in \Pi(\cc_1)$, there exists $a_n\in\cc_1$
such that $\gamma({s_n})=\Pi(a_n)$. Up to extracting a subsequence, we can assume that
$a_n\la a_\infty\in\ov{\cc}_1$. Then necessarily $\Pi(a_\infty)=\gamma(t_0)=\Pi(b')$ for some
$b'\in \ov{\cc}_2$. By properties (a) and (b), we must have
$$H(\gamma(t_0))= V_{b'}=V_{a_\infty}=\min_{c\in \llbracket a_\infty,b' \rrbracket} V_c\leq V_{c_0}.$$
This completes the proof. \endproof

\begin{proposition}
\label{identi-cc}
Let $r>0$ and $a,b\in\t_\zeta$. Then $\Pi(a)$ and $\Pi(b)$ belong to the same
connected component of $\{x\in\D_z:H(x)>r\}$ if and only if $a$ and $b$ belong
to the same connected component of $\{c\in\t_\zeta:V_c>r\}$.
\end{proposition}

\proof If $a$ and $b$ belong
to the same connected component of $\{c\in\t_\zeta:V_c>r\}$, then the line segment 
$\llbracket a,b\rrbracket$ is contained in $\{c\in\t_\zeta:V_c>r\}$, and 
$\Pi(\llbracket a,b\rrbracket)$ provides a path going from $\Pi(a)$ to $\Pi(b)$
that stays in $\{x\in\D_z:H(x)>r\}$, by property (a).

Conversely, if $a$ and $b$ belong
to different connected components of $\{c\in\t_\zeta:V_c>r\}$, then
$$\min_{c\in \llbracket a,b\rrbracket} V_c\leq r,$$
and, by Lemma \ref{cactus}, any continuous path from $\Pi(a)$
to $\Pi(b)$ must visit a point $x$ with $H(x)\leq r$. It follows that 
$\Pi(a)$ and $\Pi(b)$ belong to different
connected components of $\{x\in\D_z:H(x)>r\}$. \endproof

\noindent{\it Proof of Theorem \ref{main3}.} By Proposition \ref{identi-cc},
for every $r\geq 0$,
the projection $\Pi$ induces a one-to-one correspondence between 
connected components of $\{c\in\t_\zeta:V_c>r\}$ and
connected components of $\{x\in\D_z:H(x)>r\}$. Furthermore, 
let $\mathcal{D}$ be a connected component of $\{x\in\D_z:H(x)>r\}$,
and let $\mathcal{C}$ be the associated connected component of $\{c\in\t_\zeta:V_c>r\}$
(such that $\Pi(\cc)=\mathcal{D}$). Together with property (a) above, the fact that $\Pi$ maps the volume measure $\mathrm{vol}$
to $\mathbf{V}$ immediately shows that the boundary size $|\partial \mathcal{D}|$ can be defined
by the approximation in Theorem \ref{main3}, and that $|\partial \mathcal{D}|=|\partial\mathcal{C}|$. 
Theorem \ref{main3} is now a direct consequence of Theorem \ref{main2}. \endproof

\section{The law of components above a fixed level}
\label{sec:law-compo}

Our goal in this section is to prove Theorem \ref{main4}. To this end, we will first state and
prove a theorem about excursions ``above a fixed height'' for a snake trajectory
distributed according to $\N^{*,z}_0$. 

Let us fix $r\geq 0$. Let $\omega\in\S_0$ be chosen according to $\N_0$, or to $\N^{*,z}_0$, and consider all connected 
components of $\{a\in\t_\zeta: V_a>r\}$. If $\cc$ is one of these connected components,
we can represent $\cc$ and the labels on $\cc$ by a snake trajectory $\tilde \omega$,
which is defined as follows. First we observe that there is a unique $a_0\in\t_\zeta$
such that $a_0\in \partial \cc$ and every point of $\cc$ is a descendant of $a_0$.
Note that $V_{a_0}=r$, and that
the point $a_0$ cannot be a branching point (no branching point can have label $r$, $\N_0$ a.e. or $\N^{*,z}_0$ a.s.).
Hence we can make sense of the subtrajectory rooted at $a_0$, which we denote by
$\breve\omega$. Finally, we let $\wt\omega=\mathrm{tr}_0(\breve{\omega})$. 

We can define the boundary size $\z^*_0(\tilde\omega)$ of $\wt \omega$,
using Proposition \ref{pro-exit} (setting $\z^*_0(\tilde\omega)=Z^{(a)}_r$, where $a$ is an arbitrary point of $\cc$)
if $\omega$ is chosen according to $\N^{*,z}_0$, or the excursion theory of \cite{ALG}
if $\omega$ is chosen according to $\N_0$. We call the snake trajectories $\wt\omega$ obtained when varying 
$\cc$ among the connected components of $\{a\in\t_\zeta: V_a>r\}$ the excursions of $\omega$ above level $r$.

\begin{theorem}
\label{law-excu-abo}
Let $r>0$.
On the event $\{W^*(\omega)>r\}$, let $\wt\omega^1,\wt\omega^2,\ldots$ be the excursions of $\omega$ above level $r$, ranked in decreasing 
order of their boundary sizes. Write $\wt z_1>\wt z_2>\cdots$ for these boundary sizes. Then, under
$\N^{*,z}_0(\cdot\mid W^*>r)$, conditionally on the collection $(\wt z_i)_{i\geq 1}$,
the snake trajectories $\wt\omega^1,\wt\omega^2,\ldots$  are independent 
with respective distributions $\N^{*,\tilde z_1}_0,\N^{*,\tilde z_2}_0,\ldots$.
\end{theorem}

\noindent{\it Remark.} It is not immediately obvious that the boundary sizes $\wt z^1,\wt z^2,\ldots$ are distinct a.s.
This can however be deduced from the arguments of the proof below.

\proof
We will derive the theorem from the excursion theory of \cite{ALG}, and to this end we first need to
argue under $\N_0(\dd \omega)$. We write $\omega^1,\omega^2,\ldots$ for the excursions above $0$ ranked 
in decreasing order of their boundary sizes $\z^*_0(\omega^1),\z^*_0(\omega^2),\ldots$. 
Theorem 4 in \cite{ALG} then implies
that, under $\N_0$ and conditionally on $\z^*_0(\omega^1)=z_1,\z^*_0(\omega^2)=z_2,\ldots$, the excursions $\omega^1,\omega^2,\ldots$ are independent and the conditional distribution
of $\omega^i$ is $\N^{*,z_i}_0$.

Let $A$
be the event where exactly one excursion above $0$ hits $r$, and let $\omega^{i_0}$
be this excursion. It follows from the preceding observations that,
under $\N_0(\cdot\mid A)$,
the conditional distribution of $\omega^{i_0}$ knowing $\z^*_0(\omega^{i_0})=z$ is $\N^{*,z}_0(\cdot\mid W^*>r)$. 
Hence, if $\varphi$ is a bounded nonnegative measurable function 
on $\R_+$, and $h$ is a nonnegative measurable function on $\S_0$, we have
\begin{equation}
\label{excu-abo1}
\N_0\Big(\mathbf{1}_A\,\varphi(\z^*_0(\omega^{i_0}))\,\exp\Big(-\sum_{i=1}^\infty h(\wt\omega^{i})\Big)\Big)
= \N_0\Big(\mathbf{1}_A\,\varphi(\z^*_0(\omega^{i_0}))\,\N^{*,\z^*_0(\omega^{i_0})}_0\Big(\exp\Big(-\sum_{i=1}^\infty h(\wt\omega^i)\Big)\Big)\Big),
\end{equation}
where we use the notation $\wt\omega^1,\wt\omega^2,\ldots$ introduced in the theorem for the excursions above level $r$
(notice that this makes sense both under $\N_0$ and under $\N^{*,z}_0$, and that on the event $A$, the excursions
of $\omega$ and of $\omega^{i_0}$ above level $r$ are the same). 

We will now rewrite the left-hand side of \eqref{excu-abo1} in a different form. To this end (arguing under $\N_0(\dd \omega \mid W^*>r)$),
it is convenient to introduce all excursions of $\omega$ away from $r$: each such excursion $\ov\omega^i$, $i=1,2,\ldots$,
corresponds to one connected component of $\{a\in\t_\zeta:V_a\not =r\}$, but we exclude the connected component
containing the root $\rho$ (which may be represented by $\mathrm{tr}_r(\omega)$), and apart 
from this fact the definition of these excursions is exactly the same as that of excursions above level $r$. As previously, the excursions
$\ov\omega^i$, $i=1,2,\ldots$ are listed in decreasing order of their boundary sizes $\ov{z}^i:=\z^*_0(\ov{\omega}^i)$, $i=1,2,\ldots$. For every $i\geq 1$, we 
also let $\ov{\ve}^i$ be the sign of $\ov{\omega}^i$, so that the sequence $\wt\omega^1,\wt\omega^2,\ldots$ is obtained 
by keeping only the excursions $\ov{\omega}^i$ with $\ov{\ve}^i=1$ in the sequence $\ov\omega^i$, $i=1,2,\ldots$. 
Let $\mathcal{B}$ be the $\sigma$-field generated by $\mathrm{tr}_r(\omega)$, the sequence $(\ov{\ve}^i,\ov{z}^i)_{i=1,2\ldots}$ and the
excursions $\ov{\omega}^i$ for all $i$ such that $\ov{\ve}^i=-1$ (in other words the excursions below level $r$). 
By combining the special Markov property with \cite[Theorem 4]{ALG}, we get that conditionally on 
$\mathcal{B}$ the excursions $\wt\omega^i$, $i=1,2,\ldots$ are independent and the conditional distribution of
$\wt\omega^i$ is $\N^{*,\tilde z_i}_0$, where we write $\wt z_i=\z^*_0(\wt\omega^i)$ as in the statement of the theorem --- notice that the quantities $\wt z_i$
are $\mathcal{B}$-measurable.

The point is now that the event $A$ (and the variable $\mathbf{1}_A\,\varphi(\z^*_0(\omega^{i_0}))$) is $\mathcal{B}$-measurable.
In fact it is not hard to check that $A$ is determined by the knowledge of $\mathrm{tr}_r(\omega)$ and of the 
excursions below level $r$ (for $A$ to hold, no such excursion is allowed to contain a path that comes back to $0$ and 
then visits $r$ again). Thanks to this observation, we can rewrite the left-hand side of \eqref{excu-abo1} as
$$\N_0\Big(\mathbf{1}_A\,\varphi(\z^*_0(\omega^{i_0}))\,\prod_{i=1}^\infty \N^{*,\tilde z_i}_0(e^{-h})\Big).$$
By the same argument that led us to \eqref{excu-abo1}, this is also equal to
\begin{equation}
\label{excu-abo2}
\N_0\Big(\mathbf{1}_A\,\varphi(\z^*_0(\omega^{i_0}))\,\N^{*,\z^*_0(\omega^{i_0})}_0\Big(\prod_{i=1}^\infty \N^{*,\tilde z_i}_0(e^{-h})\Big)\Big).
\end{equation}
Notice that the law of $\z^*_0(\omega^{i_0})$ under $\N_0(\cdot\mid A)$ has a positive density with respect to Lebesgue measure
(under $\N_0$, the boundary sizes of excursions away from $0$ are the jumps of a $\phi$-CSBP under its excursion measure, see
\cite[Theorem 4]{ALG}). The equality between the quantity \eqref{excu-abo2} and the right-hand side of \eqref{excu-abo1}
for every function $\varphi$ implies that we have
\begin{equation}
\label{excu-abo3}
\N^{*,z}_0\Big(\exp\Big(-\sum_{i=1}^\infty h(\wt\omega^i)\Big)\Big)=\N^{*,z}_0\Big(\prod_{i=1}^\infty \N^{*,\tilde z_i}_0(e^{-h})\Big)
\end{equation}
for Lebesgue a.a. $z>0$. We claim that \eqref{excu-abo3} in fact holds for {\it every} $z>0$. To see this, we need a continuity argument. We restrict our attention to
functions $h$ of the type $h(\omega)=h_1(\z^*_0(\omega))h_2(\omega)$, where 
$h_1$ and $h_2$ are both (nonnegative and) bounded and continuous on $\R_+$ and $\S_0$ respectively, and there exists $\delta>0$ such that
$h_1(x)=0$ if $x\leq \delta$ and $h_2(\omega)=0$ if $W^*(\omega)\leq\delta$. Under these assumptions on $h$, one can verify that both sides of \eqref{excu-abo3} are left-continuous functions 
of $z$, which will yield our claim. Let us briefly explain this.  We write $g_1(z)$ and $g_2(z)$ for the left-hand side and 
the right-hand side of \eqref{excu-abo3} respectively. 
We use the scaling transformation 
$\theta_{z/z'}$ that
maps $\N^{*,z'}_0$ to $\N^{*,z}_0$ (see Section \ref{sna-mea}) to check that $g_i(z')\la g_i(z)$ as $z'\uparrow z$, for $i=1$ or $2$. We note that this scaling transformation maps excursions above level $r$ to excursions above level $r\sqrt{z/z'}$, and, for the function $g_2$, we observe that the collection $(\wt z_i)_{i\geq 1}$ is
the value at time $r$ of the growth-fragmentation process $\mathbf{X}$ of Theorem \ref{main2}, and we use the continuity properties of this
growth-fragmentation process (see in particular Corollary 4 in \cite{Ber2}). We omit a few details that are left to the reader.

Once we know that \eqref{excu-abo3} holds for a fixed $z>0$ and for a sufficiently large class of functions $h$, 
we obtain that the conditional distribution of the random point measure
$$\sum_{i=1}^\infty \delta_{\tilde\omega^i}$$
given $(\wt z_i)_{i\geq 1}$ is as prescribed in the statement of the theorem. This completes the proof. \endproof

\noindent{\it Proof of Theorem \ref{main4}.} 
We can derive Theorem \ref{main4} from Theorem \ref{law-excu-abo} by arguments very similar to those
of the proof of Theorem 38 in \cite{Disks} and, for this reason we only sketch the main steps of the proof.
As we already noticed in the proof of Theorem \ref{main3} in the previous section, the connected
components $\cc_1,\cc_2,\ldots$ (in the notation of Theorem \ref{main4}) are in one-to-one correspondence 
with the excursions $\wt\omega^1,\wt\omega^2,\ldots$, in such a way that the boundary size
of $\cc_i$ is equal to the boundary size $\wt z_i$ of $\wt\omega_i$. Following \cite[Section 8]{Disks},
we can associate a random compact metric space $\Theta(\wt\omega^i)$ with each excursion $\wt\omega^i$,
and we know, by Theorem \ref{law-excu-abo} and the main result of \cite{Disks}, that conditionally
on $(\wt z_i)_{i\geq 1}$, the random metric spaces $\Theta(\wt\omega^i)$, $i\geq 1$, are independent free
Brownian disks with respective perimeters $\wt z_i$, $i\geq 1$. So, all that remains is to show that,
for every $i\geq 1$, the random metric space $(\ov{\cc}_i,d_i)$ can be constructed in the way explained in
the statement of Theorem \ref{main4} and is isometric to $\Theta(\wt\omega^i)$. This is exactly similar to the 
proof of the identity (60) in \cite{Disks}, to which we refer for additional details. \hfill$\square$

\section{Complements}
\label{sec:complements}

\subsection{The cumulant function.} 
It is known \cite{Shi} that a (self-similar) growth-fragmentation process is characterized by a pair consisting 
of the self-similarity index (here $\alpha=-1/2$) and a cumulant function $\kappa$, which
is a convex function defined on $(0,\infty)$ possibly taking the values $+\infty$. The cumulant function 
$\kappa$ is given explicity in terms of the Laplace exponent $\psi$ and the L\'evy measure 
$\pi(\dd y)$ of the L\'evy process $\xi$ appearing in the Lamperti representation of the self-similar process describing the evolution
of the Eve particle,  via the formula
$$\kappa(p)=\psi(p)+ \int_{(-\infty,0)}(1-e^y)^p\,\pi(\dd y),\quad p>0.$$ 
This identity is used in \cite{BCK} to give an explicit formula for $\kappa(p)$ (see formula (33) in \cite{BCK}). We will now describe 
a different approach to the formula for $\kappa$, which is independent of the derivation of the Laplace exponent.
This suggests that another proof of Theorem \ref{main2} should be possible without the identification of
the law of the locally largest evolution, provided one knows a priori that the process $(\mathbf{Y}(r))_{r\geq 0}$ 
is a growth-fragmentation process --- note that Theorem \ref{law-excu-abo} does not provide enough information
for this.

In view of recovering the expression of $\kappa$, we observe that the negative values of
the cumulant function are given by the following formula \cite[Section 3]{BBCK}.
We consider the growth-fragmentation process $(\mathbf{Y}(r))_{r\geq 0}$ of Theorem \ref{main2} started 
at $(1,0,0,\ldots)$. For every $r\geq 0$, write $\mathbf{Y}(r)=(Y^1_r,Y^2_r,\ldots)$, and, for every $p\in\R$,
$$\|\mathbf{Y}(r)\|_{p}=\sum_{i=1}^\infty |Y^i_r|^p.$$
Then, for every $p>1/2$, the quantity  
\begin{equation}
\label{cum0}
\N^{*,1}_0\Big(\int_0^\infty \dd r\,\|\mathbf{Y}(r)\|_{p-1/2}\Big)
\end{equation}
is finite if and only if $\kappa(p)<0$, and is then equal to $-1/\kappa(p)$ \cite[Formula (16)]{BBCK}.

For every $i\geq 1$, let $\sigma^i_r$ be the volume of the $i$-th connected component of
$\{a\in \t_\zeta: V_a>r\}$ (for our purposes here the way the connected components are ordered is irrelevant). Let $q\in(-1,1)$. As a consequence of \eqref{re-root-rep-tec0}, we have, for $r>0$,
\begin{equation}
\label{cum1}
\N^*_0\Big(\sum_{i=1}^\infty \sigma^i_r (Y^i_r)^q\,e^{-\z^*_0}\Big)
= 2\,\N_0\Big(\int_{-\infty}^{-r} \dd b\,\z_b(\z_{b+r})^q\,e^{-\z_b}\Big).
\end{equation}
Let us consider first the left-hand side of \eqref{cum1}. Using
Theorem \ref{law-excu-abo} and the identity $\N^{*,z}_0(\sigma)=z^2$
in the second equality, we get
\begin{align*}
\N^*_0\Big(\sum_{i=1}^\infty \sigma^i_r (Y^i_r)^q\,e^{-\z^*_0}\Big)
&=\sqrt{\frac{3}{2\pi}} \int_0^\infty \dd z\,z^{-5/2}\,e^{-z}\,\N^{*,z}_0\Big(\sum_{i=1}^\infty \sigma^i_r (Y^i_r)^q\Big)\\
&=\sqrt{\frac{3}{2\pi}} \int_0^\infty \dd z\,z^{-5/2}\,e^{-z}\,\N^{*,z}_0\Big(\sum_{i=1}^\infty (Y^i_r)^{q+2}\Big)\\
&=\sqrt{\frac{3}{2\pi}} \int_0^\infty \dd z\,z^{q-1/2}\,e^{-z}\,\N^{*,1}_0(\| \mathbf{Y}(rz^{-1/2})\|_{q+2}),
\end{align*}
by scaling. If we integrate with respect to $\dd r$, we arrive at
\begin{equation}
\label{cum2}
\int_0^\infty \dd r \N^*_0\Big(\sum_{i=1}^\infty \sigma^i_r (Y^i_r)^q\,e^{-\z^*_0}\Big)
= \sqrt{\frac{3}{2\pi}} \,\Gamma(q+1)\,\N^{*,1}_0\Big(\int_0^\infty \dd r\,\|\mathbf{Y}(r)\|_{q+2}\Big).
\end{equation}
Consider then the right-hand side of \eqref{cum1}. Recalling formula \eqref{CSBP-tran}
and the fact that the process $(\z_{-r})_{r>0}$ is Markovian under $\N_0$
with the transition kernels of the $\phi$-CSBP, we get
\begin{align*}
\N_0\Big(\int_{-\infty}^{-r} \dd b\,\z_b(\z_{b+r})^q\,e^{-\z_b}\Big)
&=\int_{-\infty}^{-r} \dd b\,\N_0\Big((\z_{b+r})^{q+1}\times (1+r\sqrt{\frac{2}{3}})^{-3}\exp(-\z_{b+r}(1+r\sqrt{\frac{2}{3}})^{-2})\Big)\\
&=\N_0\Big(\int_{-\infty}^{0} \dd b\,(\z_{b})^{q+1}\times (1+r\sqrt{\frac{2}{3}})^{-3}\exp(-\z_{b}(1+r\sqrt{\frac{2}{3}})^{-2})\Big).
\end{align*}
Integrating with respect to $\dd r$, we find
\begin{equation}
\label{cum3}
\int_0^\infty \dd r\,\N_0\Big(\int_{-\infty}^{-r} \dd b\,\z_b(\z_{b+r})^q\,e^{-\z_b}\Big)
= \frac{1}{2\sqrt{2/3}}\,\N_0\Big(\int_{-\infty}^0 \dd b\,(\z_b)^{q} (1-e^{-\z_b})\Big).
\end{equation}
To compute the right-hand side, write $x^{q-1}=\Gamma(1-q)^{-1}\int_0^\infty \dd \lambda \,\lambda^{-q}\,e^{-\lambda x}$
(for $x>0$), and recall \eqref{Laplace-exit}, which gives $\N_0(\z_be^{-\lambda\z_b})= \lambda^{-3/2}(\lambda^{-1/2}-b\sqrt{2/3})^{-3}$
for $b<0$. It follows that
\begin{align*}
\N_0\Big(\int_{-\infty}^0 \dd b\,(\z_b)^{q} (1-e^{-\z_b})\Big)
&=\frac{1}{\Gamma(1-q)}\int_0^\infty \dd \lambda \,\lambda^{-q}\,\N_0\Big(\int_{-\infty}^0 \dd b\,\z_b\,(e^{-\lambda\z_b}-e^{-(\lambda +1)\z_b})\Big)\\
&=\frac{1}{\Gamma(1-q)}\int_0^\infty \dd \lambda \,\lambda^{-q}\int_{-\infty}^0 \dd b\,\Big(\lambda^{-3/2}(\lambda^{-1/2}-b\sqrt{2/3})^{-3}\\
&\qquad\qquad\qquad\qquad\qquad\qquad\qquad-(\lambda+1)^{-3/2}((\lambda+1)^{-1/2}-b\sqrt{2/3})^{-3}\Big)\\
&=\frac{1}{2\sqrt{2/3}\,\Gamma(1-q)}\int_0^\infty \dd \lambda \,\lambda^{-q}\,(\lambda^{-1/2}-(\lambda+1)^{-1/2}).
\end{align*}
The right-hand side is finite if and only if $-1/2<q<1/2$, and then an elementary calculation gives
$$\int_0^\infty \dd \lambda \,\lambda^{-q}\,(\lambda^{-1/2}-(\lambda+1)^{-1/2}) = - \frac{1}{\sqrt{\pi}}\,
\Gamma(1-q)\Gamma(q-\frac{1}{2}).$$
Coming back to \eqref{cum3}, we see that
$$\int_0^\infty \dd r\,\N_0\Big(\int_{-\infty}^{-r} \dd b\,\z_b(\z_{b+r})^q\,e^{-\z_b}\Big)
= -\frac{3}{8\sqrt{\pi}}\Gamma(q-\frac{1}{2}).$$
Combining this equality with \eqref{cum1} and \eqref{cum2} leads to
$$\N^{*,1}_0\Big(\int_0^\infty \dd r\,\|\mathbf{Y}(r)\|_{q+2}\Big)
= -\sqrt{\frac{3}{8}}\,\frac{\Gamma(q-\frac{1}{2})}{\Gamma(q+1)}.$$
Replacing $q$ by $p=q+5/2$, we finally obtain that, for $2<p<3$, the quantity \eqref{cum0} is finite, and
$$\kappa(p)= \sqrt{\frac{8}{3}}\,\frac{\Gamma(p-\frac{3}{2})}{\Gamma(p-3)},$$
which is in agreement with formula (33) in \cite{BCK} --- note that the value of 
$\kappa$ in the latter formula should be multiplied by the factor $\sqrt{3/2\pi}$ 
that appears in the formula for $\psi$ in Theorem \ref{main1}. 

Finally, an argument of analytic continuation shows that the preceding formula
for $\kappa(p)$ holds for every $p>3/2$, whereas $\kappa(p)=+\infty$ for $p\in(0,3/2]$. 
The function $p\mapsto \kappa(p)$ is (finite and) convex on $(3/2,\infty)$, and vanishes at 
$p=2$ and $p=3$ (with the notation of \cite{BBCK}, we have $\omega_-=2$ and $\omega_+=3$). 

\subsection{A growth-fragmentation process in the Brownian plane}

In this section, we consider the random pointed metric space $(\mathcal{P}_\infty,\dd_\infty)$ called the Brownian plane, which has been introduced and studied in \cite{Plane}.
The space $\mathcal{P}_\infty$ has a distinguished point $\rho_\infty$, and, for every $r>0$, we may define the boundary sizes of
the connected components of $\{x\in\mathcal{P}_\infty:\dd(\rho_\infty,x)>r\}$, via the same approximation as used above in Section \ref{sec:bdry-size}: To see that this definition makes sense, one may argue that there exists a
coupling of the Brownian plane and the Brownian map such that small balls centered at the distinguished point in the two spaces are isometric \cite{Plane}, 
then rely on Proposition \ref{approx-exit} to treat the case when $r$ is small enough, and finally 
use the scale invariance of the Brownian plane.
Notice that there is exactly one unbounded component, whose
boundary is also the boundary of the so-called hull of radius $r$ in $\mathcal{P}_\infty$ (see in particular \cite{CLG}).

We will relate this collection of boundary sizes to the growth-fragmentation process of Theorem \ref{main1} subject 
to a special conditioning. Precisely, we consider this growth-fragmentation process starting from $0$ and conditioned to
have indefinite growth (see \cite[Section 4.2]{BBCK}). Let us briefly describe this process, referring to \cite{BBCK} for more details. We start with one Eve particle, whose mass process
$(\wh X_t)_{t\geq 0}$ evolves as the process $X$ of Theorem \ref{main1}
conditioned to start from $0$ and to stay positive for all times. To be specific, the process $\wh X$ is a self-similar Markov process
with index $1/2$, which can be obtained via the Lamperti representation from a L\'evy process $\wh\xi$ with no positive jumps and Laplace exponent
\begin{equation}
\label{psi-condit}
\wh\psi(\lambda):=\kappa(3+\lambda)= \sqrt{\frac{8}{3}}\,\frac{\Gamma(\lambda+\frac{3}{2})}{\Gamma(\lambda)},\quad \lambda>0.
\end{equation}
See \cite[Lemma 2.1]{BBCK} for the fact that the function $\wh\psi(\lambda)$ corresponds to the Laplace exponent of 
a L\'evy process without positive jumps.
Then, as previously, each jump time $t$ of $\wh X$ corresponds to the birth of a new particle (child of the Eve particle) with mass $|\Delta\wh X_t|$,
but the masses of these new particles evolve independently according to the distribution of the process $X$, and similarly for the children 
of these particles, and so on. We emphasize that only the mass process of the Eve particle evolves according to a different Markov process $\wh X$, while 
the masses of its children, grandchildren, etc., evolve according to the law of the process $X$. As previously, we write $\wh{\mathbf{X}}(r)$ for the
collection of masses of all particles present at time $r$.

\begin{theorem}
\label{growth-frag-plane}
As a process indexed by the variable $r> 0$, the collection of the boundary sizes of
all connected components of $\{x\in\mathcal{P}_\infty:\dd_\infty(\rho_\infty,x)>r\}$
is distributed as the process $(\wh{\mathbf{X}}(r))_{r> 0}$.
\end{theorem}

\proof We first explain that the role of the Eve particle (for the process $\wh{\mathbf{X}}$) is played by the evolution of
the unbounded component of $\{x\in\mathcal{P}_\infty:\dd_\infty(\rho_\infty,x)>r\}$. Let $\wh Z_r$ be the boundary size of this component, with the
convention that $\wh Z_0=0$.
The distribution of the process $(\wh Z_r)_{r\geq 0}$ is described in \cite[Proposition 1.2]{CLG}. From this
description, using also the arguments of \cite[Section 4.4]{CLG0}, one gets that 
$\wh Z_r$ can be written as
$$\wh Z_r = U^\uparrow_{\eta_r}$$
where $(U^\uparrow_t)_{t\geq 0}$ is the L\'evy process with no positive jumps and Laplace exponent $\phi(\lambda)=\sqrt{8/3}\,\lambda^{3/2}$
conditioned to start from $0$ and to stay positive for all times $t>0$ --- see \cite[Chapter VII]{Ber} for a rigorous definition of this
process --- and $\eta_r$ is the time change
$$\eta_r=\inf\{t\geq 0: \int_0^t \frac{\dd s}{U^\uparrow_s} >r\}.$$
It follows from this representation that $(\wh Z_r)_{r\geq 0}$ is a self-similar Markov process with index $1/2$ with values 
in $[0,\infty)$. We will now verify that the Laplace exponent of the L\'evy process 
arising in the Lamperti representation of this self-similar Markov process is equal to $\wh\psi(\lambda)$, which
will imply that the process $(\wh Z_t)_{t\geq 0}$ has the same distribution as the mass process of the Eve particle 
in the description of the process $(\wh{\mathbf{X}}(t))_{t\geq 0}$. We slightly abuse notation by introducing, for every $x\geq 0$, a
 probability measure $\P_x$ under which the Markov process $\wh Z$ starts from $x$. By self-similarity, for every $a>0$, the law of $(a^{-2}\wh Z_{at})_{t\geq 0}$
under ${\P}_x$ coincides with the law of $(\wh{Z}_t)_{t\geq 0}$ under $\P_{a^{-2}x}$. We recall from \cite[Proposition 1.2]{CLG} that, for every
$r>0$, $\wh Z_r$ follows (under $\P_0$) a Gamma distribution with parameter $3/2$ and mean $r^2$.

Let $q\in(-\frac{3}{2},-\frac{1}{2})$. Then,
\begin{equation}
\label{gfp1}
\E_0\Big[\int_1^\infty \wh Z_t^q\,\dd t\Big]= \int_1^\infty t^{2q}\,\E_0[\wh Z_1^q]\,\dd t= - \frac{\E_0[\wh Z_1^q]}{2q +1}= -\left(\frac{2}{3}\right)^q \frac{1}{2q+1} \frac{\Gamma(q+\frac{3}{2})}
{\Gamma(\frac{3}{2})}.
\end{equation}
We may compute the left-hand side of \eqref{gfp1} in a different manner by applying the Markov property at time $1$. We get
\begin{equation}
\label{gfp2}
\E_0\Big[\int_1^\infty \wh Z_t^q\,\dd t\Big]= \E_0\Big[ \E_{\wh Z_1}\Big[\int_0^\infty \wh Z_t^q\,\dd t\Big]\Big] = \E_0\Big[ \wh Z_1^{q+1/2}\Big] \times 
\E_{1}\Big[\int_0^\infty \wh Z_t^q\,\dd t\Big]
\end{equation}
using the self-similarity of $\wh Z$. Then $\E_0[ \wh Z_1^{q+1/2}]= (2/3)^{q+1/2}\Gamma(q+2)/\Gamma(3/2)$
and, on the other hand, if $\wh \xi$ denotes the L\'evy process (started from $0$) arising in the Lamperti representation of the
self-similar process $\wh Z$, we have
$$\E_{1}\Big[\int_0^\infty \wh Z_t^q\,\dd t\Big] = \E\Big[ \int_0^\infty e^{(q+\frac{1}{2})\hat \xi(t)}\,\dd t\Big].$$
The quantity in the right-hand side must be finite, which implies that $\E[e^{(q+\frac{1}{2})\hat \xi(t)}]<\infty$ for every $t>0$, and
$\E[e^{(q+\frac{1}{2})\hat \xi(t)}]=\exp(t\psi^*(q+\frac{1}{2}))$ with $\psi^*(q+\frac{1}{2})<0$. From  \eqref{gfp1} and \eqref{gfp2}, we get
$$\frac{1}{\psi^*(q+\frac{1}{2})} 
= \frac{1}{2} \left(\frac{2}{3}\right)^{-1/2}\, \frac{\Gamma(q+\frac{1}{2})}{\Gamma(q+2)}.$$
Finally, we find that, for $-1<q<0$, we have
$$\psi^*(q)= \sqrt{\frac{8}{3}}\, \frac{\Gamma(q+\frac{3}{2})}{\Gamma(q)}=\wh \psi(q).$$
An argument of analytic continuation now allows us to obtain that the Laplace exponent 
of the L\'evy process $\wh \xi$ is equal to $\wh\psi(\lambda)$ as desired. 

Once we have identified $(\wh Z_t)_{t\geq 0}$ as the mass process of the Eve particle 
in the description of $(\wh{\mathbf{X}}(t))_{t\geq 0}$, the remaining steps of the proof
are very similar to those of the proof of Theorem \ref{main2}, and we will only
sketch the main ingredients. We first recall the relevant features of the construction of the Brownian plane $(\mathcal{P}_\infty,\dd_\infty)$ 
which is developed in \cite[Section 3.2]{CLG}, to which we refer for further details. The random metric space $\mathcal{P}_\infty$ is obtained as
a quotient space of a (non-compact) random tree $\t_\infty$, which itself is constructed by
grafting a Poisson collection of (compact) $\R$-trees to an infinite spine isometric to $[0,\infty)$. The point $0$ of the spine corresponds to the
distinguished point $\rho_\infty$ of $\mathcal{P}_\infty$. Furthermore, every point $a$ of $\t_\infty$
is assigned a nonnegative label $\Lambda_a$, and this label coincides with $\dd_\infty(\rho_\infty,x)$, if $x$
is the point of $\mathcal{P}_\infty$ corresponding to $a$.
Then, as in the proof of Theorem \ref{main3}, it is not hard to check that connected components of $\{x\in\mathcal{P}_\infty: \dd_\infty(\rho_\infty,x)>r\}$
are in one-to-one correspondence with connected components of $\{a\in\t_\infty: \Lambda_a>r\}$, for every fixed $r>0$.

For every $a\in\t_\infty$, let  $\llbracket a,\infty\llbracket$ stand for the range of the unique geodesic path from $a$
to $\infty$ in $\t_\infty$, and set $\un{\Lambda}_a=\min\{\Lambda_b: b\in\llbracket a,\infty\llbracket\}$. If $\cc$ is a (necessarily bounded)
connected component of $\{a\in\t_\infty: \Lambda_a-\un\Lambda_a>0\}$, then both $\cc$
and the labels $(\Lambda_a)_{a\in\cc}$ can be represented by a snake trajectory $\omega_\cc$, 
in a way very similar to what we did for $\check\cc_r$ in Section \ref{excu-mini}.

\begin{proposition}
\label{excu-infinite-tree}
Setting $\inf\{\Lambda_a:a\in\cc\}=r$ yields a one-to-one 
correspondence between
connected components $\cc$ of $\{a\in\t_\infty: \Lambda_a-\un\Lambda_a>0\}$ and jump times $r$ of the process $(\wh Z_t)_{t\geq 0}$. Let $r_1,r_2,\ldots$ be an enumeration of these
jump times, which is measurable with respect to the $\sigma$-field generated by $(\wh Z_r)_{r\geq 0}$, and for every $i=1,2,\ldots$, let 
$\cc_i$ be the connected component  associated with $r_i$.  Then, conditionally 
on the process $(\wh Z_r)_{r\geq 0}$, the snake trajectories $\omega_{\cc_i}$, $i=1,2,\ldots$, are independent, and
the conditional distribution of $\omega_{\cc_i}$ is $\N^{*,|\Delta \wh Z_{r_i}|}_0$. 
\end{proposition}

Proposition \ref{excu-infinite-tree} is obviously an analog of Theorem \ref{mainALG}, and can be derived 
from the latter result using the relations between the labeled tree $(\t_\infty,(\Lambda_a)_{a\in\t_\infty})$
and the pair $(\t_\zeta,(V_a)_{a\in\t_\zeta})$ under $\N_0$ (compare the decomposition of the Brownian snake
at the minimum found in \cite[Theorem 2.1]{CLG} with the construction of $(\t_\infty,(\Lambda_a)_{a\in\t_\infty})$
in the same reference). 

Given Proposition \ref{excu-infinite-tree}, the end of the proof of Theorem \ref{growth-frag-plane} follows the same
general pattern as that of Theorem \ref{main2} in Section \ref{gro-frag}, and we leave the details to the
reader. \endproof

\noindent{\it Remark.} We could have used Corollary 2 of \cite{CabCha} to identify the Laplace exponent of the L\'evy process $\wh\xi$ in
the preceding proof. This would still have required some calculations, and for this reason we preferred the 
more direct approach presented above.
 
\subsection{Local times}

In this section, we argue under $\N_0(\dd\omega)$, but similar results hold under $\N^*_0(\dd \omega)$.
Recall from the introduction the definition of the local times $(\mathcal{L}_x,x\in\R)$
of the process $(V_a)_{a\in\t_\zeta}$.
In this section, we 
relate the values of $\mathcal{L}_x$ for $x>0$ to the growth-fragmentation process of Theorem \ref{main1} or equivalently
to the connected components of $\{a\in\t_\zeta:V_a>x\}$. 

We fix $h>0$. It is convenient to introduce the  ``local time exit process'' $(\mathscr{X}^h_r)_{r\geq 0}$, which roughly 
speaking measures for every $r\geq 0$ the ``quantity'' of paths $W_s$ that have hit level $h$ and accumulated 
a local time equal to $r$ at level $h$. The precise definition of this process fits in the general framework of exit measures
\cite[Chapter V]{Zurich} and we refer to the introduction of \cite{ALG} for more details (there only the case $h=0$
is considered, but the extension to the case $h>0$ is straightforward). Note that
$\mathscr{X}^h_0=\z_h$ is just the usual exit measure from $(-\infty,r)$, which can be defined by formula \eqref{formu-exit}.
Furthermore, under $\N_0(\cdot \mid W^*>h)$, the process $(\mathscr{X}^h_r)_{r\geq 0}$ is a $\phi$-CSBP 
started from $\z_h$  --- 
see again the introduction of \cite{ALG}. Of course, on the event $\{W^*<h\}$, the process $(\mathscr{X}^h_r)_{r\geq 0}$ is identically zero.

Let $\mathcal{C}^h_1,\mathcal{C}^h_2,\ldots$ be the connected components of $\{a\in\t_\zeta:V_a>h\}$, ranked 
in decreasing order of their boundary sizes $|\partial\cc^h_1|>|\partial\cc^h_2|>\cdots$. For every
$i=1,2,\ldots$, let $H_i$ be the height of $\cc^h_i$, defined by
$$H_i:=\sup_{a\in\cc^h_i} V_a - h.$$

\begin{proposition}
\label{local-approx}
We have $\N_0$ a.e.
\begin{equation}
\label{ident-local}
\mathcal{L}_h= \int_0^\infty \dd r \, \mathscr{X}^h_r.
\end{equation}
Moreover, $\N_0$ a.e.,
\begin{equation}
\label{conv10}
\delta^{3/2}\#\{i: |\partial\cc^h_i|>\delta\}\build{\la}_{\delta\to 0}^{}\frac{1}{\sqrt{6\pi}}\,\mathcal{L}_h
\end{equation}
and 
\begin{equation}
\label{conv20}
\delta^3\#\{i: H_i >\delta\}\build{\la}_{\delta\to 0}^{} \mathbf{c}\,\mathcal{L}_h
\end{equation}
where
$\mathbf{c}= \frac{3}{2}\pi^{-3/2}\Gamma(1/3)^3\Gamma(7/6)^3$.
\end{proposition}

\noindent{\it Remark.} The proposition also holds for $h=0$, but the proof of \eqref{ident-local} requires a slightly different argument
in that case. We omit the details.

\proof 
The convergence \eqref{conv20} is already established in \cite{Cactus}, in the slightly weaker form
of convergence in measure (note that ``upcrossings'' from $h$ to $h+\delta$, as defined in \cite{Cactus},
exactly correspond to connected components of $\{a\in\t_\zeta:V_a>h\}$ with height greater than $\delta$). 
We will use this fact to prove the identity \eqref{ident-local}. To simplify notation, we set
$$\mathcal{L}^*_h=\int_0^\infty \dd r \, \mathscr{X}^h_r.$$
As in the proof of Theorem \ref{law-excu-abo}, we consider all excursions away from $h$.
It follows from \cite[Proposition 3, Theorem 4]{ALG} (and an application of the special Markov property) that these excursions
are in one-to-one correspondence with the jumps of the process $(\mathscr{X}^h_r)_{r\geq 0}$, and that conditionally
on the latter process, they are independent and the conditional distribution of the excursion
corresponding to the jump $\Delta\mathscr{X}^h_r$ is 
$$\frac{1}{2}\N^{*,\Delta\mathscr{X}^h_r}_0 + \frac{1}{2}\check\N^{*,\Delta\mathscr{X}^h_r}_0,$$
where we use the notation $\check\N^{*,z}_0$ for the push forward of $\N^{*,z}_0$ under the symmetry $\omega\mapsto -\omega$. 
We recall that the connected components of $\{a\in\t_\zeta:V_a>h\}$ are in one-to-one correspondence 
with the excursions above level $h$, in such a way that the boundary size of a component is equal to
the corresponding jump of  $(\mathscr{X}^h_r)_{r\geq 0}$. 

Write $U$ for the (stopped) L\'evy process obtained from the $\phi$-CSBP $\mathscr{X}^h$ by the Lamperti transformation (note that
$U$ is stopped upon hitting $0$ and that $U_0=\z_h$). Notice that the hitting time of $0$ by $U$ is $\mathcal{L}^*_h$.
Since the jumps
of $(\mathscr{X}^h_r,)_{r\geq 0}$ are also the jumps of $U$, we obtain the identity in distribution
$$\Big(\mathcal{L}^*_h, \sum_{i=1}^\infty \delta_{|\partial \cc^h_i|}\Big) \build{=}_{}^{\rm(d)} \Big(\mathcal{L}^*_h, \sum_{i=1}^\infty \mathbf{1}_{\{\epsilon_i=1\}}\delta_{\Delta U_{r_i}}\Big)$$
where $r_1,r_2,\ldots$ are the jump times of $U$, and $\epsilon_1,\epsilon_2,\ldots$ is a sequence of independent Bernoulli variables 
with parameter $1/2$, which is independent of $U$. Since the L\'evy measure of $U$ is
$$\pi(\dd z)=\sqrt{\frac{3}{2\pi}}\,z^{-5/2}\,\dd z,$$
so that $\pi((\delta,\infty))=\sqrt{2/3\pi}\,\delta^{-3/2}$, it easily follows that, $\N_0$ a.e.,
$$\delta^{3/2}\#\{i: \epsilon_i=1\hbox{ and }\Delta U_{r_i}>\delta\}\build{\la}_{\delta\to 0}^{}\frac{1}{\sqrt{6\pi}}\,\mathcal{L}^*_h.$$
This gives the convergence \eqref{conv10}, except that we have not yet verified that $\mathcal{L}_h=\mathcal{L}^*_h$.

To this end, using again the conditional distribution of the excursions away from $h$ given the 
process $(\mathscr{X}^h_r)_{r\geq 0}$, we observe that we have also
$$\Big(\mathcal{L}^*_h, \sum_{i=1}^\infty \delta_{H_i}\Big) \build{=}_{}^{\rm(d)} \Big(\mathcal{L}^*_h, \sum_{i=1}^\infty \mathbf{1}_{\{\epsilon_i=1\}}\delta_{
\sqrt{\Delta U_{r_i}}\,M_i}\Big)$$
where $M_1,M_2,\ldots$ is a sequence of independent random variables distributed according to the law of $W^*$
under $\N^{*,1}_0$, which is also independent of $(U,(\epsilon_i)_{i\geq 1})$. Now observe that, if $z$ is chosen according to $\pi(\dd z)$ and $M$
according to the law of $W^*$
under $\N^{*,1}_0$, $\sqrt{z}\,M$ is distributed according to the ``law'' of $W^*$ under $\N^*_0$, which satisfies
$$\N^*_0(W^*>\delta)= 2\mathbf{c}\,\delta^{-3}$$
by \cite[Lemma 25]{ALG}. It follows that, $\N_0$ a.e.,
$$\delta^3\#\{i: H_i >\delta\}\build{\la}_{\delta\to 0}^{}\mathbf{c}\,\mathcal{L}^*_h.$$
By comparing this convergence with \cite[Theorem 6]{Cactus}, we get that $\mathcal{L}_h=\mathcal{L}^*_h$, which completes the proof.  \endproof

\section*{Appendix}


\smallskip
\noindent{\it Proof of Lemma \ref{Lamperti-excu}}. It is convenient to write $\mathbf{N}$ 
for the distribution of $(\z_{-t})_{t\geq 0}$ under $\N_0$ (we agree that $\z_0=0$). Then $\mathbf{N}$
is a $\sigma$-finite measure on the Skorokhod space $\mathbb{D}([0,\infty),\R)$. 
For $\ve >0$, let 
$$\sum_{i\in I_\ve} \delta_{w^\ve_i}(\dd w)$$
be a Poisson point measure on $\mathbb{D}([0,\infty),\R)$ with intensity $\ve\mathbf{N}$. 
As already noticed in Section \ref{exit-mea}, we can construct a 
$\phi$-CSBP started from $\ve$
by setting, for $t>0$,
$$Y^\ve_t=\sum_{i\in I_\ve} w^\ve_i(t)$$
and $Y^\ve_0=\ve$. 
Set $T^\ve_0:=\inf\{t\geq 0: Y^\ve_t=0\}$.
The classical Lamperti transformation \cite{Lam0,CLU} allows us to relate
$Y^\ve$ to another process $X^\ve$ distributed as a stable L\'evy process 
with no negative jumps and index $3/2$, started from $\ve$ and stopped at 
the first time when it hits $0$, via the formula
$$X^\ve_t=Y^\ve_{\gamma^\ve_t}$$
where $\gamma^\ve_t:=\inf\{s\geq 0:\int_0^s Y^\ve_u\dd u>t\}$ if $t<\int_0^{\infty} Y^\ve_u\dd u$ and
$\gamma^\ve_t=T^\ve_0$ otherwise. 

Let us fix $\delta>0$ and assume from now on that $\ve\in(0,\delta)$. 
Let $B_\ve$ stand for the event $\{\sup_{t\geq 0}Y^\ve_t\geq \delta\}=\{\sup_{t\geq 0}X^\ve_t\geq \delta\}$. By the solution of the two-sided exit problem already used in the proof of Lemma \ref{law-Lz},
we have 
$$\P(B_\ve)= 1- \sqrt{\frac{\delta-\ve}{\delta}} = \frac{\ve}{2\delta} +O(\ve^2)\;,\quad\hbox{as }\ve\to 0.$$
On the other hand, let $A_\ve\subset B_\ve$ be the event where there is exactly one $i\in I_\ve$
such that $\sup_{t\geq 0} w^\ve_i(t)\geq \delta$. By Lemma \ref{law-Lz},
$$\P(A_\ve) = \frac{\ve}{2\delta}\,\exp\Big(- \frac{\ve}{2\delta}\Big)= \frac{\ve}{2\delta} +O(\ve^2)\;,\quad\hbox{as }\ve\to 0.$$

If $F$ is a bounded measurable function on the space $\mathbb{D}([0,\infty),\R)$, we deduce from the last two displays that
\begin{equation}
\label{app1}
\E[ F((X^\ve_t)_{t\geq 0})\mid B_\ve] = \E[ F((X^\ve_t)_{t\geq 0})\mid A_\ve] + O(\ve)\;,\quad\hbox{as }\ve\to 0.
\end{equation}

We can associate with $X^\ve$ the process ``reflected above the minimum'' defined by
$$\tilde X^\ve_t:= X^\ve_t - \inf\{X^\ve_s:0\leq s\leq t\}.$$
We have obviously $0\leq X^\ve_t-\tilde X^\ve_t\leq \ve$ for every $t\geq 0$. If $\tilde B_\ve\subset B_\ve$
stands for the event $\{\sup_{t\geq 0}\tilde X^\ve_t\geq \delta\}$, it is easily checked that
$\P(B_\ve\backslash \tilde B_\ve)=O(\ve^2)$, so that we can replace $B_\ve$ by $\tilde B_\ve$
in \eqref{app1}. Furthermore, on the event $\tilde B_\ve$ we can introduce the first excursion
of $\tilde X^\ve$ away from $0$ that hits $\delta$ and denote this excursion by $(\mathcal{X}^\ve_t)_{t\geq 0}$. 
Notice that the distribution of $(\mathcal{X}^\ve_t)_{t\geq 0}$ is $\mathbf{n}_\delta(\dd e):=\mathbf{n}(\dd e\mid\sup\{e(t):t\geq 0\}\geq \delta)$.
Let $\dd_{\mathrm{Sk}}$ be a distance on $\mathbb{D}([0,\infty),\R)$ that induces the Skorokhod topology. It is a
simple matter to verify that, for every $\alpha >0$,
$$\P(\dd_{\mathrm{Sk}}((\tilde X^\ve_t)_{t\geq 0},(\mathcal{X}^\ve_t)_{t\geq 0})>\alpha\mid \tilde B_\ve) \build{\la}_{\ve \to 0}^{} 0.$$
Assume from now on that $F$ is (bounded and) Lipschitz with respect to $\dd_{\mathrm{Sk}}$. We deduce from 
\eqref{app1} (with $B_\ve$ replaced by $\tilde B_\ve$) and the last display that
\begin{equation}
\label{app2}
\E[ F((X^\ve_t)_{t\geq 0})\mid A_\ve]  - \E[ F((\mathcal{X}^\ve_t)_{t\geq 0})\mid \tilde B_\ve] \build{\la}_{\ve \to 0}^{} 0.
\end{equation}
Note that $\E[ F((\mathcal{X}^\ve_t)_{t\geq 0})\mid \tilde B_\ve] $ does not depend on $\ve$ and is equal to $\int\mathbf{n}_\delta(\dd e)\,F(e)$. 

On the other hand, conditionally on the event $A_\ve$ there is a unique index $i_0\in I_\ve$ such that
$\sup_{t\geq 0} w^\ve_{i_0}(t)\geq \delta$, and the distribution of 
$w^\ve_{i_0}$ is $\mathbf{N}_\delta(\dd w):=\mathbf{N}(\dd w\mid \sup_{t\geq 0} w(t)\geq \delta)$. We then set $\mathcal{Y}^\ve_t=w^\ve_{i_0}(t)$, and
$$\eta^\ve_t=\inf\{s\geq 0: \int_0^s \dd u\,\mathcal{Y}^\ve_u>t\}$$
if $t< \int_0^\infty \dd u\,\mathcal{Y}^\ve_u$, and $\eta^\ve_t=\inf\{s\geq 0: \mathcal{Y}^\ve_s=0\}$ otherwise. 
Observing that the conditional distribution of $Y^\ve-\mathcal{Y}^\ve$ given $A_\ve$ is dominated 
by the law of a $\Phi$-continuous state branching process started from $\ve$, one verifies that, for every $\alpha>0$,
$$\P(\dd_{\mathrm{Sk}}((Y^\ve_{\gamma^\ve_t})_{t\geq 0},(\mathcal{Y}^\ve_{\eta^\ve_t})_{t\geq 0})>\alpha\mid A_\ve) \build{\la}_{\ve \to 0}^{} 0.$$
(Here we omit a few details that are left to the reader.) Recalling that $Y^\ve_{\gamma^\ve_t}=X^\ve_t$ and using \eqref{app2}, we get
$$\E[ F((\mathcal{Y}^\ve_{\eta^\ve_t})_{t\geq 0})\mid A_\ve]  - \int\mathbf{n}_\delta(\dd e)\,F(e)\build{\la}_{\ve \to 0}^{} 0.$$
Since the conditional distribution of $(\mathcal{Y}^\ve_{\eta^\ve_t})_{t\geq 0}$ given $A_\ve$ is $\mathbf{N}_\delta$ (independently of 
$\ve$), using the equalities $\mathbf{N}(\{w:\sup_{t\geq 0}w(t)\geq \delta\})=\frac{1}{2\delta}=\mathbf{n}(\{e:\sup_{t\geq 0}e(t)\geq \delta\})$
(the first one by Lemma \ref{law-Lz} and the second one as an easy consequence of the two-sided exit problem), we arrive at the
result of the lemma. \endproof

\noindent{\it Proof of Proposition \ref{approx-exit}}. {\it First step.} Recall that, for every $t>0$,  $\g_t$ denotes the $\sigma$-field 
on $\S_0$ generated by the mapping $\omega\mapsto \mathrm{tr}_{-t}(\omega)$ and completed by the 
collection of all $\N_0$-negligible sets. We also define $\g_0$ as the $\sigma$-field generated 
by the $\N_0$-negligible sets. For every $\eta>0$, the process $(\z_{-t})_{t\geq \eta}$  is Markov with respect to the filtration $(\g_{t})_{t\geq \eta}$ 
under the probability measure $\N_0^{[\eta]}:=\N_0(\cdot \mid W_*\leq -\eta)$. By the Feller property
of the semigroup, the strong Markov property holds even for stopping times of 
the filtration $(\g_{t+})_{t\geq \eta}$.

We fix two reals $\eta\in(0,1)$ and $M>1$. Let $\ve\in (0,\eta)$. From the proof of Proposition 34 in the appendix of \cite{Disks},
we have, for every $r\leq -\eta$,
$$\N_0\Big((\z^\ve_r - \z_r)^2\Big)\leq 4\ve^2.$$
We note that \cite{Disks} deals with the quantity $\tilde \z^\ve_r$ defined in Remark (ii) after
Proposition \ref{approx-exit}, rather than with $\z^\ve_r$, but as explained in this remark, this makes no difference for
a fixed value of $r$. 
Furthermore, \cite{Disks} gives the latter bound only for ``truncated versions'' of $\tilde\z^\ve_r$ and $\z_r$, but an application of
Fatou's lemma then yields the preceding display. 

Let $\delta\in(0,1)$. By Markov's inequality, for $r\leq -\eta$,
$$\N_0(|\z^\ve_r-\z_r|>\delta) \leq \delta^{-2}\times 4\ve^2.$$
We apply this to $r=-j\,\ve^{3/2}$ for all integers $j$ such that $\eta\leq j\,\ve^{3/2}\leq M+1$. It follows that
\begin{align}
\label{app10}
&\N_0\Big(|\z^\ve_{-j\ve^{3/2}}-\z_{-j\ve^{3/2}}|>\frac{\delta}{2},\hbox{ for some }j\hbox{ s.t. } \eta\leq j\,\ve^{3/2}\leq M+1\Big)\\
&\qquad\leq 16(M+1)\delta^{-2}\,\ve^{1/2}.\nonumber
\end{align}

Fix a real $K>0$, and consider the random time
$$S:=\inf\{t\geq \eta: \z_{-t}<K\hbox{ and } |\z^\ve_{-t}-\z_{-t}| >\delta\}.$$
Note that $S$ is a stopping time of the filtration $(\g_{t+})_{t\geq \eta}$ (because 
both processes $(\z^\ve_{-t})_{t\geq \eta}$ and $(\z_{-t})_{t\geq \eta}$ have c\`adl\`ag paths and are
adapted to the filtration $(\g_{t})_{t\geq \eta})$. On the event $\{S<\infty\}$, we have 
$ |\z^\ve_{-S}-\z_{-S}| \geq\delta$ and $\z_{-S}\leq K$. 

Our goal is now to bound $\N_0(S\leq M)$. To this end, we will use \eqref{app10}. On the event
$\{S<\infty\}$ write $[-S]_\ve$ for the greatest number of the form $-j\ve^{3/2}$ in the interval $(-\infty,-S)$. Then,
$$\{S\leq M\} = \{S\leq M,|\z^\ve_{[-S]_\ve}-\z_{[-S]_\ve}|>\delta/2\} \cup \{S\leq M,|\z^\ve_{[-S]_\ve}-\z_{[-S]_\ve}|\leq\delta/2\}.$$
By \eqref{app10}, the $\N_0$-measure of the first set in the right-hand side is bounded above by $c_1\ve^{1/2}$ for some constant $c_1$ depending on $M$ and $\delta$. On the other hand,
recalling that $ |\z^\ve_{-S}-\z_{-S}| \geq\delta$ on $\{S<\infty\}$, we obtain that the second set is contained in
$$\{S\leq M,\,|\z_{[-S]_\ve}-\z_{-S}|\geq\delta/4\}\cup \{S\leq M,\,|\z^\ve_{[-S]_\ve}-\z^\ve_{-S}|\geq\delta/4\}.$$
Using the strong Markov property of $(\z_{-t})_{t\geq \eta}$ at time $S$, the bound $\z_{-S}\leq K$ on $\{S<\infty\}$,
and the fact that a $\phi$-CSBP can be written 
as a time change of a L\'evy process,  it is easy to verify that
\begin{equation}
\label{app11}
\N_0(S\leq M,\,|\z_{[-S]_\ve}-\z_{-S}|\geq\delta/4)\leq c_2\ve^{3/2}
\end{equation}
for some constant $c_2$ depending on $\delta$ and $K$. 

In the second and the third step below, we will get similar estimates 
for the $\N_0$-measure of (appropriate subsets of) the event $\{S\leq M,\,|\z^\ve_{[-S]_\ve}-\z^\ve_{-S}|\geq\delta/4\}$.
We will explain in the fourth step how the proof of the proposition is completed by combining all these estimates.

\medskip
\noindent{\it Second step.} We first study the quantity
$$\N_0(S\leq M,\,\z^\ve_{[-S]_\ve}-\z^\ve_{-S}\geq\delta/4).$$
From our definitions, on the event $\{S<\infty\}$, the quantity 
$\z^\ve_{[-S]_\ve}-\z^\ve_{-S}$ is bounded above by
$$F_\ve:=\ve^{-2}\,\int_0^\sigma \dd s\,\mathbf{1}_{\{T_{-S}(W_s)<\infty,T_{-S-\ve^{3/2}}(W_s)=\infty\}}\,\mathbf{1}_{\{\wh W_s<[-S]_\ve+\ve\}}.$$
For every integer $n\geq 1$, write $[-S]_{(n)}$ for the greatest number of the form $-j2^{-n}$ in $(-\infty,-S)$, and set, still on the event $\{S<\infty\}$,
$$F_{\ve,n}:=\ve^{-2}\,\int_0^\sigma \dd s\,\mathbf{1}_{\{T_{[-S]_{(n)}}(W_s)<\infty,T_{[-S]_{(n)}-\ve^{3/2}}(W_s)=\infty\}}\,\mathbf{1}_{\{\wh W_s<[-S]_{(n)}+\ve\}}.$$
Observing that 
$\mathbf{1}_{\{T_{-S}(W_s)<\infty\}}=\lim_{n\to\infty}\mathbf{1}_{\{T_{[-S]_{(n)}}(W_s)<\infty\}}$ and using Fatou's lemma, we have
$$\N_0(\mathbf{1}_{\{S<\infty\}}\,F_{\ve})\leq \liminf_{n\to\infty} \N_0(\mathbf{1}_{\{S<\infty\}}\,F_{\ve,n}).$$
Then,
\begin{align*}
&\N_0(\mathbf{1}_{\{S<\infty\}}\,F_{\ve,n})\\
&\quad=\ve^{-2}\sum_{k=1}^\infty \N_0\Big(\mathbf{1}_{\{(k-1)2^{-n}\leq S<k2^{-n}\}}
\int_0^\sigma \dd s\,\mathbf{1}_{\{T_{-k2^{-n}}(W_s)<\infty,T_{-k2^{-n}-\ve^{3/2}}(W_s)=\infty\}}\,\mathbf{1}_{\{\wh W_s<-k2^{-n}+\ve\}}\Big).
\end{align*}
We can apply the special Markov property (Proposition \ref{SMP}) to each term of the sum in the right-hand side.
Note that the variable $\mathbf{1}_{\{(k-1)2^{-n}\leq S<k2^{-n}\}}$ is measurable with respect to $\g_{k2^{-n}}$, whereas 
the subsequent integral is a function of the snake trajectories $\omega_i$ introduced in Proposition \ref{SMP} when $r=-k2^{-n}$). 
We obtain
$$
\N_0(\mathbf{1}_{\{S<\infty\}}\,F_{\ve,n})=\ve^{-2}\sum_{k=1}^\infty \N_0\Big(\mathbf{1}_{\{(k-1)2^{-n}\leq S<k2^{-n}\}}\,
\z_{k2^{-n}}\Big) 
\N_0\Big(\int_0^\sigma \dd s\,\mathbf{1}_{\{T_{-\ve^{3/2}}(W_s)=\infty\}}\,\mathbf{1}_{\{\wh W_s<\ve\}}\Big).
$$
By the first-moment formula for the Brownian snake \cite[Proposition 4.2]{Zurich}, we have 
$$\N_0\Big(\int_0^\sigma \dd s\,\mathbf{1}_{\{T_{-\ve^{3/2}}(W_s)=\infty\}}\,\mathbf{1}_{\{\wh W_s<\ve\}}\Big)=
\E_0\Big[\int_0^{\mathfrak{t}_{-\ve^{3/2}}} \dd t\,\mathbf{1}_{\{B_t<\ve\}}\Big] \leq c_3\ve^{5/2},$$
where $(B_t)_{t\geq 0}$ is a standard linear Brownian motion starting from $x$ under the probability measure $\P_x$, $\mathfrak{t}_{r}=\inf\{t\geq 0:B_t=r\}$ 
for every $r\in\R$, and $c_3$ is a constant.
We conclude that $\N_0(\mathbf{1}_{\{S<\infty\}}\,F_{\ve,n})\leq c_3\ve^{1/2}\,\N_0(\mathbf{1}_{\{S<\infty\}}\,\z_{[S]_{(n)}})\leq c_3\ve^{1/2}$, and
the same bound holds for $\N_0(\mathbf{1}_{\{S<\infty\}}\,F_{\ve})$. Finally, Markov's inequality gives
$$
\N_0(S\leq M,\,\z^\ve_{[-S]_\ve}-\z^\ve_{-S}\geq\delta/4)\leq \N_0(S<\infty,\,F_\ve\geq\delta/4)
\leq \frac{4}{\delta}\,c_3\ve^{1/2}.
$$

\noindent{\it Third step.} We now consider the event
$\{S\leq M,\,\z^\ve_{-S}-\z^\ve_{[-S]_\ve}\geq\delta/4)\}$.
We observe that, if 
$S<\infty$,
$$\z^\ve_{-S}-\z^\ve_{[-S]_\ve}\leq \ve^{-2}\int_0^\sigma \dd s\,\mathbf{1}_{\{T_{-S}(W_s)=\infty,\,\wh W_s\in [[-S]_\ve +\ve,-S+\ve)\}}.$$
Notice that $T_{-S}(W_s)=\infty$ implies $T_{[-S]_\ve}(W_s)=\infty$ and that $-S+\ve\leq [-S]_\ve +\ve +\ve^{3/2}$.
Hence, on the event where $S\leq M$ and $\z^\ve_{-S}-\z^\ve_{[-S]_\ve}\geq\delta/4$, we can find a real $r\in[-M-1,-\eta]$
of the form $r=j\ve^{3/2}$ with $j\in \mathbb{Z}$, such that
$$\ve^{-2}\int_0^\sigma \dd s\,\mathbf{1}_{\{T_{r}(W_s)=\infty,\,\wh W_s\in [r +\ve,r+\ve+\ve^{3/2}]\}}\geq \frac{\delta}{4}.$$

Let us fix $r\in[-M-1,-\eta]$ in the following lines, and bound the probability of the event in the last display. 
From the first-moment formula for the Brownian snake, we have, with the same notation as in the second step,
$$
\N_0\Big(\ve^{-2}\int_0^\sigma \dd s\,\mathbf{1}_{\{T_{r}(W_s)=\infty,\,\wh W_s\in [r +\ve,r+\ve+\ve^{3/2}]\}}\Big)
=\ve^{-2}\E_0\Big[\int_0^{\mathfrak{t}_r} \dd t\,\mathbf{1}_{\{B_t\in [r+\ve,r+\ve+\ve^{3/2}]\}}\Big]\leq c_4\ve^{1/2},
$$
with some constant $c_4$. To get a better estimate, we use higher moments, but to this end we need to perform a
suitable truncation. We fix $A>0$, and we observe that, for every integer $k\geq 1$, for any nonnegative measurable function $f$ on $\R$, we have
$$
\N_0\Big(\Big( \int_0^\sigma \dd s \,\mathbf{1}_{\{T_r(W_s)=\infty,\tau_A(W_s)=\infty\}} f(\wh W_s)\Big)^k\Big)
\leq C_{k,A,M}\, \Big(\sup_{x\in[r,A]} \E_x\Big[\int_0^{\mathfrak{t}_r\wedge \mathfrak{t}_A} \dd t\,f(B_t)\Big]\Big)^k,
$$
where $C_{k,A,M}$ is a constant depending only on $k$, $A$ and $M$. The bound in the previous display can be derived in
a straightforward way from the $k$-th moment formula for the Brownian snake \cite[Proposition IV.2]{Zurich}. We omit the details.
We apply this bound with $f=\mathbf{1}_{ [r+\ve,r+\ve+\ve^{3/2}]}$, and we arrive at the estimate
$$\N_0\Big(\Big(\ve^{-2}\int_0^\sigma \dd s \,\mathbf{1}_{\{T_r(W_s)=\infty,\tau_A(W_s)=\infty\}} \mathbf{1}_{\{\wh W_s\in [r+\ve,r+\ve+\ve^{3/2}]\}}\Big)^k\Big)
\leq C_{k,A,M}\,(c_5\ve^{1/2})^k,$$
with a constant $c_5$ depending on $A$ and $M$. From Markov's inequality, we then get
\begin{align*}
&\N_0\Big(W^*<A,\;\ve^{-2}\int_0^\sigma \dd s\,\mathbf{1}_{\{T_{r}(W_s)=\infty,\,\wh W_s\in [r +\ve,r+\ve+\ve^{3/2}]\}}\geq \frac{\delta}{4}\Big)\\
&\quad\leq \left(\frac{\delta}{4}\right)^{-k}\N_0\Big(\Big(\ve^{-2}\int_0^\sigma \dd s \,\mathbf{1}_{\{T_r(W_s)=\infty,\tau_A(W_s)=\infty\}}
\mathbf{1}_{\{\wh W_s\in [r+\ve,r+\ve+\ve^{3/2}]\}}\Big)^k\Big)\\
&\quad \leq C_{k,A,M}\,\left(\frac{\delta}{4}\right)^{-k}\,(c_5\ve^{1/2})^k.
\end{align*}
We take $k=4$ and sum the preceding estimate over possible values of $r=-j\ve^{3/2}$ in $[-M-1,-\eta]$, and we arrive at the estimate
$$\N_0(W^*<A,\,S\leq M,\,\z^\ve_{-S}-\z^\ve_{[-S]_\ve}\geq\delta/4)\leq c_6\,\ve^{1/2}$$
with a constant $c_6$ depending on $A,M$ and $\delta$. 

\medskip
\noindent{\it Fourth step.} We deduce from the second and third steps that we have
\begin{equation}
\label{app12}
\N_0(W^*<A,\,S\leq M,\,|\z^\ve_{[-S]_\ve}-\z^\ve_{-S}|\geq\delta/4)|\leq c_7\,\ve^{1/2},
\end{equation}
with a constant $c_7$ depending on $\delta, A, M, K$. Combining \eqref{app11} and \eqref{app12} and recalling
the considerations of the end of the first step, we arrive at the bound 
\begin{equation}
\label{app13}
\N_0(W^*<A,\,S\leq M)\leq (c_1+c_2+c_7)\,\ve^{1/2}.
\end{equation}
Let us write $S=S^{(\ve)}$ to recall the dependence on $\ve$. Let $n_0$ be the first integer such that $(n_0)^{-3}<\eta$. The bound 
\eqref{app13} gives
$$\sum_{n=n_0}^\infty \N_0(W^*<A,\,S^{(n^{-3})}\leq M)<\infty.$$
Hence, $\N_0$ a.e. on the event $W^*<A$, we have $S^{(n^{-3})}> M$ for all large enough $n$. This means that, $\N_0$ a.e. on the event 
where $\sup\{\z_{-t}:t> 0\}<K$ and $W^*<A$,
we have for all large enough $n$,
$$\sup_{\eta\leq u\leq M} |\z^{(n^{-3})}_{-u}-\z_{-u}|\leq \delta.$$
Since  $\delta$, $K$ and $A$ are arbitrary, we obtain that, $\N_0$ a.e.,
$$\lim_{n\to\infty} \Big(\sup_{\eta\leq u\leq M} |\z^{(n^{-3})}_{-u}-\z_{-u}|\Big)=0.$$
We can replace $\eta\leq u\leq M$ by $\eta\leq u<\infty$ since $\z^{(\ve)}_{-u}-\z_{-u}=0$ for 
$u>-W_*+\ve$.
The statement of the proposition then follows by a monotonicity argument. \hfill$\square$

\smallskip\noindent
{\bf Acknowledgement.} We are indebted to Nicolas Curien for a number of very stimulating discussions.

\end{document}